%% file: main.tex
\pdfoutput=1
\documentclass[reprint]{JASA}

\input{settings/packages.tex}
\input{settings/typesetting.tex}

\input{settings/mathcommands.tex}

\input{settings/mathsymbols.tex}

\input{resources/glossary.tex}
\input{resources/abbreviations.tex}

\input{settings/supmat-commands.tex}

\begin{document}

\input{head/titlepage.tex}

\glsresetall{}

\input{head/abstract.tex}
\glsresetall{}

\maketitle

\input{content/introduction.tex}

\input{content/approach.tex}

\input{content/validation.tex}

\input{content/discussion.tex}

\input{content/conclusion.tex}

\input{tail/acknowledgments.tex}

\bibliography{resources/references.bib}

\supplements
\input{content/sup-introduction.tex}
\clearpage
\input{content/sup-approach.tex}

\clearpage
\input{content/sup-validation.tex}

\end{document}

%% file: settings/packages.tex
\usepackage{csquotes} %
\usepackage{textcomp}  %
\usepackage{siunitx}  %
\DeclareSIUnit[number-unit-product = {}]{\percent}{\char 37}

\usepackage[caption=false]{subfig}
\newcommand{\phantomsubfloat}[1]{%
    {%
        \captionsetup[subfloat]{farskip=0pt,captionskip=0pt}%
        \captionsetup[subfigure]{labelformat=empty}%
        \subfloat[][]{#1}%
    }%
}
\newcommand{\phantomsubfloatprevspace}{\vspace{-.33333333333\baselineskip}}

\usepackage{mathtools}  %
\usepackage{cases}

\usepackage{multirow}

\usepackage{glossaries-extra}
\setabbreviationstyle[acronym]{long-short}  %
\glsdisablehyper

\usepackage[capitalize]{cleveref} %
\crefname{table}{Table}{Tables} %
\crefname{figure}{Fig.\@}{Fig.\@} %
\Crefname{figure}{Figure}{Figures}
\crefname{section}{Sec.\@}{Secs.\@}
\Crefname{section}{Section}{Sections} %
\crefname{equation}{Eq.\@}{Eqs.\@}
\Crefname{equation}{Equation}{Equations} %
\crefname{algorithm}{Algorithm}{Algorithms} %
\Crefname{algorithm}{Algorithm}{Algorithms} %

%% file: settings/typesetting.tex
\newcommand*{\ndim}[1]{\mbox{#1-D}}  %

%% file: settings/mathcommands.tex
\DeclarePairedDelimiter{\parens}{\lparen}{\rparen}

\DeclarePairedDelimiter{\bracks}{\lbrack}{\rbrack}

\DeclarePairedDelimiter{\abs}{\lvert}{\rvert}

\DeclarePairedDelimiterXPP\onenorm[1]{}{\lVert}{\rVert}{_{1}}{#1}
\DeclarePairedDelimiterXPP\twonorm[1]{}{\lVert}{\rVert}{_{2}}{#1}
\DeclarePairedDelimiterXPP\fronorm[1]{}{\lVert}{\rVert}{_{\text{F}}}{#1}
\DeclarePairedDelimiterXPP\infnorm[1]{}{\lVert}{\rVert}{_{\infty}}{#1}
\DeclarePairedDelimiterXPP\pnorm[1]{}{\lVert}{\rVert}{_{p}}{#1}
\DeclarePairedDelimiterXPP\tvnorm[1]{}{\lVert}{\rVert}{_{\text{TV}}}{#1}
\DeclarePairedDelimiterXPP\maxnorm[1]{}{\lVert}{\rVert}{_{\text{max}}}{#1}
\DeclarePairedDelimiterXPP\nuclearnorm[1]{}{\lVert}{\rVert}{_{\ast}}{#1}

\DeclareMathOperator{\sinc}{sinc}

\newcommand{\dirac}{\delta}

\DeclareMathOperator*{\conv}{\ast}

\ifLuaTeX
    \newcommand{\mat}[1]{\symbfit{#1}}
    \renewcommand{\vec}[1]{\symbfit{#1}}
\else
    \newcommand{\mat}[1]{\bm{#1}}
    \renewcommand{\vec}[1]{\bm{#1}}
\fi

\newcommand{\inv}{^{-1}}

\DeclareMathOperator{\determinant}{det}

\newcommand{\crossprod}{\times}

\newcommand{\dd}{\mathop{}\!\mathrm{d}}

\newcommand{\set}[1]{\mathbb{#1}}

\newcommand{\integernumbers}{\set{Z}}

\newcommand{\realnumbers}{\set{R}}

\newcommand{\ellipsis}{\dots}
\let\setroster\braces

\newcommand{\estimate}[1]{\hat{#1}}

\DeclarePairedDelimiter{\tuple}{\lparen}{\rparen}

\DeclarePairedDelimiter\floor{\lfloor}{\rfloor}

%% file: settings/mathsymbols.tex
\newcommand{\soundspeedletter}{c}
\newcommand{\soundspeed}{\soundspeedletter}
\newcommand{\soundspeedmean}{\soundspeed}

\newcommand{\densityletter}{\rho}
\newcommand{\density}{\densityletter}
\newcommand{\densitymean}{\density}

\newcommand{\wavelength}{\lambda}

\newcommand{\velpotletter}{\psi}
\newcommand{\velpot}{\velpotletter}

\newcommand{\pressureletter}{p}
\newcommand{\pressure}{\pressureletter}

\newcommand{\surfsepletter}{\xi}

\newcommand{\waveformsurface}{v_{\surface}}
\newcommand{\waveformsurfacedistrib}{\surfsepletter}

\newcommand{\sourcespecifier}[1]{#1'}
\newcommand{\fieldpointletter}{r}
\newcommand{\fieldpoint}{\vec{\fieldpointletter}}
\newcommand{\fieldpointsource}{\sourcespecifier{\fieldpoint}}
\newcommand{\surfnormalletter}{n}
\newcommand{\surfnormal}{\vec{\surfnormalletter}}

\newcommand{\timevartletter}{t}
\newcommand{\timevar}{\timevartletter}

\newcommand{\sourceimagespecifier}[1]{#1}

\newcommand{\surfnormalsourceimage}{\sourcespecifier{\sourceimagespecifier{\surfnormal}}}

\newcommand{\ddsurf}{\dd \sigma}
\newcommand{\surfaceletter}{S}
\newcommand{\surface}{\surfaceletter}
\newcommand{\volumeletter}{V}
\newcommand{\volume}{\volumeletter}

\newcommand{\sirsymb}{h}
\newcommand{\sir}{\sirsymb}

\newcommand{\sirbctermsymb}{\beta}
\newcommand{\sirbcterm}{\sirbctermsymb}
\newcommand{\sirsurface}{\sir_{\surface}}
\newcommand{\sirsurfaceintegrand}{g}

\newcommand{\sirdelay}{\tau}

\newcommand{\nurbsknotvectoru}{U}
\newcommand{\nurbsknotvectorv}{V}
\newcommand{\nurbscoordu}{u}
\newcommand{\nurbscoordv}{v}
\newcommand{\nurbsdegreeu}{p}
\newcommand{\nurbsdegreev}{q}
\newcommand{\nurbsbasis}{b}
\newcommand{\nurbsindexu}{i}
\newcommand{\nurbsindexv}{j}
\newcommand{\nurbsintervala}{0}  %
\newcommand{\nurbsintervalb}{1}  %
\newcommand{\nurbsbasisnbu}{n}  %
\newcommand{\nurbsbasisnbv}{m}  %
\newcommand{\nurbsknotsnbu}{r}  %
\newcommand{\nurbsknotsnbv}{s}  %
\newcommand{\nurbsctrlpoint}{\vec{p}}
\newcommand{\nurbsweight}{w}

\newcommand{\nurbssurfacemap}{\vec{s}}
\newcommand{\nurbssurface}{\surface}

\newcommand{\surfacemap}{\vec{s}}
\newcommand{\surfaceparamspace}{\hat{S}}
\newcommand{\surfaceparamvec}{\vec{u}}
\newcommand{\surfacemapjac}{\mat{J}_{\surfacemap}}
\newcommand{\quadpointindex}{q}
\newcommand{\quadweightsymb}{w}
\newcommand{\quadpointnb}{n_{\quadpointindex}}
\newcommand{\surfacequadcoordvec}{\surfaceparamvec_{\quadpointindex}}
\newcommand{\surfacequadcoordu}{u}
\newcommand{\surfacequadcoordv}{v}
\newcommand{\surfacequadpoint}{\fieldpoint_{\quadpointindex}}
\newcommand{\surfacequadweight}{\quadweightsymb_{\quadpointindex}}
\newcommand{\surfacequadjacdet}{j_{\quadpointindex}}
\newcommand{\surfacequadnormal}{\surfnormal_{\quadpointindex}}

\newcommand{\sirgenericweight}{\alpha}

\newcommand{\sirgenericsignal}{y}
\newcommand{\sirgenericsignalvec}{\vec{\sirgenericsignal}}

\newcommand{\sirgenericexcitation}{v}
\newcommand{\sirgenericexcitationcoeffs}{c}
\newcommand{\sirexcitationtimeindex}{k}
\newcommand{\sirexcitationtimenb}{K}
\newcommand{\sirdeltat}{T}

\newcommand{\interpbasisletter}{\varphi}
\newcommand{\interpbasis}{\interpbasisletter_{\text{int}}}
\newcommand{\interpbasisgeneral}{\interpbasisletter}
\newcommand{\interpbasisconvinv}{\interpbasisgeneral\inv}

\newcommand{\interpbspline}{\beta}
\newcommand{\interpbsplinedegree}{n}
\newcommand{\interpbsplinepolepairnb}{m}
\newcommand{\interpbsplinepole}{z}
\newcommand{\interpbsplinepoleindex}{i}

\newcommand{\sirsurfacegeneral}{\hat{\sir}_{\surface}}

\newcommand{\apodletter}{a}
\newcommand{\sirarrayelemindex}{i}
\newcommand{\sirarrayelemnb}{n_{e}}
\newcommand{\sirarrayelemwaveform}{\sirgenericexcitation_{\sirarrayelemindex}}
\newcommand{\sirarrayelemcoeffs}{\sirgenericexcitationcoeffs_{\sirarrayelemindex}}
\newcommand{\sirarrayelembsir}{\hat{\sir}_{\sirarrayelemindex}}
\newcommand{\sirarrayelemdelay}{\sirdelay_{\sirarrayelemindex}}
\newcommand{\sirarrayelemapod}{\apodletter_{\sirarrayelemindex}}

\newcommand{\sirexpdiracnb}{n}
\newcommand{\sirexpdiracindex}{i}
\newcommand{\sirexpdiracweight}{\apodletter}
\newcommand{\sirexpdiracdelay}{\sirdelay}
\newcommand{\sirexpwaveform}{\sirgenericexcitation}
\newcommand{\sirexpsignal}{\sirgenericsignal}

\DeclareMathOperator{\lognormaldistsymb}{Lognormal}
\newcommand{\lognormaldist}{\lognormaldistsymb}

%% file: resources/glossary.tex
\newglossaryentry{python}{%
    name={Python},
    description={interpreted, high-level, general-purpose programming language}
}
\newglossaryentry{matlab}{%
    name={MATLAB},
    description={multi-paradigm numerical computing environment and proprietary programming language developed by MathWorks}
}
\newglossaryentry{tensorflow}{%
   name={TensorFlow},
   description={An end-to-end open source machine learning platform}
}
\newglossaryentry{tensorrt}{%
   name={TensorRT},
   description={NVIDIA TensorRT™ is a platform for high-performance deep learning inference.}
}
\newglossaryentry{mu-law}{%
    name={\textmu-law},
    description={Comanding algorithm}
}
\newglossaryentry{ge}{%
    name={General Electric Company},
    description={TODO}
}

\newglossaryentry{ge-healthcare}{%
    name={GE Healthcare},
    description={TODO}
}

\newglossaryentry{ge9ld}{%
    name={9L-D},
    description={GE Linear transducer}
}
\newglossaryentry{bmode}{%
    name={B-mode},
    description={Brightness mode}
}
\newglossaryentry{verasonics}{%
    name={Verasonics},
    description={TODO}
}
\newglossaryentry{vantage256}{%
    name={Vantage~256},
    description={TODO}
}

\newglossaryentry{cirs}{%
    name={CIRS},
    description={Multi-purpose multi-tissue ultrasound phantom}
}

\newglossaryentry{cirs-040gse}{%
    name={\gls{cirs} model~040GSE},
    description={Multi-purpose multi-tissue ultrasound phantom}
}
\newglossaryentry{cirs-054gs}{%
    name={\gls{cirs} model~054GS},
    description={General purpose ultrasound phantom}
}
\newglossaryentry{pyus}{%
    name={PyUS},
    description={TODO}
}
\newglossaryentry{fieldii}{%
    name={Field~II},
    description={TODO}
}
\newglossaryentry{pivlab}{%
    name={PIVlab},
    description={TODO}
}
\newglossaryentry{famusii}{%
    name={FAMUS~II},
    description={TODO}
}
\newglossaryentry{kwave}{%
    name={k-Wave},
    description={TODO}
}
\newglossaryentry{x-ray}{%
    name={X-ray},
    description={X-radiation or Röntgen radiation}
}
\newglossaryentry{nvidia}{%
    name={NVIDIA},
    description={American technology company incorporated in Delaware and based in Santa Clara, California. It designs graphics processing units (GPUs) for the gaming and professional markets, as well as system on a chip units (SoCs) for the mobile computing and automotive market.}
}
\newglossaryentry{1080ti}{
    name={NVIDIA GeForce GTX 1080~Ti},
    description={Specific GPU from NVIDIA (Pascal architecture)}
}
\newglossaryentry{mx150}{
    name={NVIDIA GeForce MX~150},
    description={Specific GPU from NVIDIA (Pascal architecture)}
}
\newglossaryentry{titanv}{
    name={NVIDIA TITAN~V},
    description={Specific GPU from NVIDIA (Volta architecture)}
}
\newglossaryentry{v100}{
    name={NVIDIA Tesla V100},
    description={Specific GPU from NVIDIA (Volta architecture)}
}
\newglossaryentry{unet}{
    name={U-Net},
    description={TODO}
}
\newglossaryentry{resnet}{
    name={ResNet},
    description={Typical neural network architecture}
}
\newglossaryentry{he}{
    name={He},
    description={TODO He init}
}
\newglossaryentry{glorot}{
    name={Glorot},
    description={TODO Gloro init}
}
\newglossaryentry{bspline}{
    name={B-spline},
    description={TODO}
}
\newglossaryentry{bernstein}{
    name={Bernstein},
    description={TODO}
}
\newglossaryentry{bezier}{
    name={Bézier},
    description={TODO}
}
\newglossaryentry{fourier}{
    name={Fourier},
    description={TODO}
}
\newglossaryentry{keys}{
    name={Keys},  %
    description={TODO}
}
\newglossaryentry{sinc}{
    name={sinc},
    description={TODO}
}
\newglossaryentry{z-transform}{
    name={Z-transform},
    description={TODO}
}
\newglossaryentry{simpson}{
    name={Simpson},
    description={TODO}
}
\newglossaryentry{jacobian}{
    name={Jacobian},
    description={TODO}
}
\newglossaryentry{etc}{
    name={etc.},
    description={et cetera},
}
\newglossaryentry{ie}{
    name={i.e.}, description={id est}, %
}
\newglossaryentry{eg}{
    name={e.g.}, description={exempli gratia}, %
}
\newglossaryentry{etal}{
    name={\textit{et al.}}, description={et alii}, %
}
\newglossaryentry{ibid}{
    name={\textit{ibid.}},
    description={ibidem}, %
}
\newglossaryentry{vs}{
    name={vs.}, description={TODO}, %
}
\newglossaryentry{wrt}{
    name={with respect to}, description={TODO}, %
}
\newglossaryentry{esp}{
    name={esp.}, description={TODO}, %
}
\newglossaryentry{invitro}{
    name={\textit{in vitro}}, description={TODO}, %
}
\newglossaryentry{invivo}{
    name={\textit{in vivo}}, description={TODO}, %
}
\newglossaryentry{exvivo}{
    name={\textit{ex vivo}}, description={TODO}, %
}
\newglossaryentry{insitu}{
    name={\textit{in situ}}, description={TODO}, %
}
\newglossaryentry{insilico}{
    name={\textit{in silico}}, description={TODO}, %
}
\newglossaryentry{interalia}{
    name={\textit{inter alia}}, description={TODO}, %
}
\newglossaryentry{intoto}{
    name={\textit{in toto}}, description={TODO}, %
}
\newglossaryentry{apriori}{
    name={\textit{a priori}}, description={TODO}, %
}
\newglossaryentry{aposteriori}{
    name={\textit{a posteriori}}, description={TODO}, %
}
\newglossaryentry{alas}{
    name={\textit{alas}}, description={TODO}, %
}
\newglossaryentry{defacto}{
    name={de facto}, description={TODO}, %
}
\newglossaryentry{viceversa}{
    name={\textit{vice versa}}, description={TODO}, %
}

%% file: resources/abbreviations.tex
\newacronym{nn}{NN}{neural network}
\newacronym{dnn}{DNN}{deep neural network}
\newacronym{cnn}{CNN}{convolutional neural network}
\newacronym{relu}{ReLU}{rectified linear unit}
\newacronym{conv-layer}{CL}{convolutional layer}
\newacronym{std-conv-block}{FCB}{fully convolutional block}
\newacronym{res-conv-block}{RCB}{residual convolutional block}
\newacronym{mse}{MSE}{mean squared error}
\newacronym{mae}{MAE}{mean absolute error}
\newacronym{mslae}{MSLAE}{mean signed logarithmic absolute error}
\newacronym{mmuae}{MMUAE}{mean \textmu-law absolute error}
\newacronym{roi}{ROI}{region of interest}
\newacronym{voi}{VOI}{volume of interest}
\newacronym{fps}{FPS}{frames per second}
\newacronym{snr}{SNR}{signal-to-noise ratio}
\newacronym{cnr}{CNR}{contrast-to-noise ratio}
\newacronym{gcnr}{GCNR}{generalized contrast-to-noise ratio}
\newacronym{fwhm}{FWHM}{full width at half maximum}
\newacronym{ssim}{SSIM}{structural similarity}
\newacronym{psnr}{PSNR}{peak signal-to-noise ratio}
\newacronym{dra}{DRA}{dynamic range alteration}
\newacronym{epe}{EPE}{endpoint error}
\newacronym{mepe}{MEPE}{mean \glsxtrlong*{epe}}
\newacronym{repe}{REPE}{relative \glsxtrlong*{epe}}
\newacronym{mrepe}{MREPE}{mean \glsxtrlong*{repe}}
\newacronym{rve}{RVE}{ratio of valid estimates}
\newacronym{mri}{MRI}{magnetic resonance imaging}
\newacronym{ct}{CT}{computed tomography}
\newacronym{dc}{DC}{direct current}
\newacronym{hdr}{HDR}{high dynamic range}
\newacronym{nurbs}{NURBS}{non-uniform rational B-spline}
\newacronym{fft}{FFT}{fast Fourier transform}
\newacronym{ifft}{IFFT}{inverse fast Fourier transform}
\newacronym{fir}{FIR}{finite impulse response}
\newacronym{iir}{IIR}{infinite impulse response}
\newacronym{moms}{MOMS}{maximal-order-minimal-support}
\newacronym{o-moms}{O-MOMS}{optimal-maximal-order-minimal-support}
\newacronym{cad}{CAD}{computer-aided design}
\newacronym{iga}{IGA}{isogeometric analysis}
\newacronym{bp}{BP}{backprojection}
\newacronym{fbp}{FBP}{filtered backprojection}
\newacronym{fista}{FISTA}{fast iterative shrinkage-thresholding algorithm}
\newacronym{sr}{SR}{sparse regularization}
\newacronym{us}{US}{ultrasound}
\newacronym{das}{DAS}{delay-and-sum}
\newacronym{trf}{TRF}{tissue reflectivity function}
\newacronym{pw}{PW}{plane wave}
\newacronym{sa}{SA}{synthetic aperture}
\newacronym{dw}{DW}{diverging wave}
\newacronym{cpwc}{CPWC}{coherent plane wave compounding}
\newacronym{rf}{RF}{radio frequency}
\newacronym{prf}{PRF}{pulse repetition frequency}
\newacronym{tgc}{TGC}{time gain compensation}
\newacronym{lq}{LQ}{low-quality}
\newacronym{hq}{HQ}{high-quality}
\newacronym{uq}{UQ}{ultra-high-quality}
\newacronym{picmus}{PICMUS}{plane-wave imaging challenge in medical ultrasound}
\newacronym{sir}{SIR}{spatial impulse response}
\newacronym{psf}{PSF}{point spread function}
\newacronym{iq}{IQ}{in-phase quadrature}
\newacronym{gl}{GL}{grating lobe}
\newacronym{sl}{SL}{side lobe}
\newacronym{ew}{EW}{edge wave}
\newacronym{tx}{Tx}{transmit}
\newacronym{rx}{Rx}{receive}
\newacronym{mv}{MV}{minimum variance}
\newacronym{mi}{MI}{mechanical index}
\newacronym{ussr}{USSR}{ultrasound sparse regularization framework}
\newacronym{piv}{PIV}{particle image velocimetry}
\newacronym{pdf}{PDF}{probability density function}
\newacronym{cdf}{CDF}{cumulative distribution function}
\newacronym{wss}{WSS}{wide-sense stationary}
\newacronym{acf}{ACF}{autocovariance function}
\newacronym{pcc}{PCC}{Pearson correlation coefficient}
\newacronym{zncc}{ZNCC}{zero-normalized cross-correlation}
\newacronym{epfl}{EPFL}{École polytechnique fédérale de Lausanne}
\newacronym{lts5}{LTS5}{Signal Processing Laboratory 5}
\newacronym{chuv}{CHUV}{University Hospital Center}
\newacronym{unil}{UNIL}{University of Lausanne}
\newacronym{cibm}{CIBM}{Center for Biomedical Imaging}
\newacronym{gpu}{GPU}{graphics processing unit}
\newacronym{cpu}{CPU}{central processing unit}
\newacronym{ieee}{IEEE}{Institute of Electrical and Electronics Engineers}
\newacronym{ius}{IUS}{International Ultrasonic Symposium}
\newacronym{sa-vfi}{SA-VFI}{synthetic aperture vector flow imaging}

%% file: head/titlepage.tex
\newcommand*{\affiliationlts}{%
    Signal Processing Laboratory 5 (LTS5),
    École polytechnique fédérale de Lausanne (EPFL),
    1015 Lausanne,
    Switzerland%
}
\newcommand*{\affiliationchuv}{%
    Department of Radiology,
    University Hospital Center (CHUV),
    and
    University of Lausanne (UNIL),
    1011 Lausanne,
    Switzerland%
}
\newcommand*{\affiliationcibm}{%
    CIBM Center for Biomedical Imaging,
    1015 Lausanne,
    Switzerland
}

\title[Perdios \textit{et al.}/Preprint]{A spline-based spatial impulse response simulator}
\author{Dimitris Perdios}
\email{dimitris.perdios@epfl.ch}
\author{Florian Martinez}
\author{Marcel Arditi}
\affiliation{\affiliationlts}
\author{Jean-Philippe Thiran}
\email{jean-philippe.thiran@epfl.ch}
\altaffiliation{Also at: \affiliationchuv.}
\altaffiliation{Also at: \affiliationcibm.}
\affiliation{\affiliationlts}

\date{\today}

%% file: head/abstract.tex
\begin{abstract}

The \gls{sir} method is a well-known approach to calculate transient acoustic
fields of arbitrary-shape transducers.
It involves the evaluation of a time-dependent
surface integral.
Although analytic expressions of the \gls{sir} exist for some geometries,
numerical methods based on the discretization of transducer surfaces
have become the standard.
The proposed method consists of representing the transducer
as a \gls{nurbs} surface,
and decomposing it into smooth \gls{bezier} patches
onto which quadrature rules can be deployed.
The evaluation of the \gls{sir} can then be expressed in \gls{bspline} bases,
resulting in a sum of shifted-and-weighted basis functions.
Field signals are eventually obtained by a convolution
of the basis coefficients,
derived from the excitation waveform,
and the basis \gls{sir}.
The use of \gls{nurbs} enables exact representations of common transducer
elements.
High-order Gaussian quadrature rules enable high accuracy
with few quadrature points.
High-order \gls{bspline} bases are ideally suited to exploit efficiently
the bandlimited property of excitation waveforms.
Numerical experiments demonstrate that the proposed approach enables
sampling the \gls{sir} at low sampling rates,
as required by the excitation waveform,
without introducing additional errors on simulated field signals.

\end{abstract}

%% file: content/introduction.tex
\section{Introduction}%
\label{sec:introduction}

The need for fast and accurate \gls{us} simulation tools is as great as ever.
Such tools have been and continue to be extensively used
for the design and characterization
of \gls{us} transducers~\cite{Kino_BOOK_1987,Hunt_TBME_1983}.
The use of end-to-end \gls{us} scanner simulators is also increasing
for the training and evaluation of physicians
and sonographers~\cite{Ziv_SIH_2006,Lewiss_JUM_2014},
as they eliminate the need for volunteers or patients and can provide
on-demand exposure to specific care and diagnosis scenarios.
Simulation tools are also crucial for the development, assessment, and validation
of image analysis and image reconstruction methods.
This includes, for example,
the analysis and characterization of speckle
patterns~\cite{Wagner_TSU_1983,Foster_UI_1983},
the optimization of application-specific acquisition
sequences~\cite{JensenJonas_UFFC_2016},
the fine-tuning and benchmarking of image reconstruction
parameters and methods~\cite{Rindal_UFFC_2019},
or the development of image quality
metrics~\cite{Rodriguez-Molares_UFFC_2020}.
The recent breakthrough of deep learning-based methods in medical image
analysis~\cite{Greenspan_TMI_2016},
medical image reconstruction~\cite{Wang_ACCESS_2016,Wang_TMI_2018},
and ultrasound imaging~\cite{vanSloun_JPROC_2020}
comes with a critical need for large-scale datasets
to feed the data-intensive algorithms involved.
In this context,
the generation of synthetic data is of great interest~\cite{Frangi_TMI_2018},
with the potential to generate (infinitely) large,
highly diverse, and unbiased datasets.

With the myriad of \gls{us} applications in which simulations can be
leveraged,
it is no wonder that many \gls{us} simulators are available.%
\footnote{A fairly exhaustive list of available software can be found
on the \gls{kwave} website: \url{http://www.k-wave.org/acousticsoftware.php}.}
Among them,
and in the context of linear acoustics,
\gls{fieldii}~\cite{Jensen_MBEC_1996,Jensen_UFFC_1992}
is considered a reference,
thanks to its ability to calculate acoustic fields of arbitrarily shaped,
excited, and apodized transducers,
as well as simulating the acquisition of ultrasound pulse-echo signals
in the presence of (weak) tissue inhomogeneities~\cite{Jensen_JASA_1991}.
It is based on the \gls{sir} method developed by
\citet{Tupholme_MATHEMATIKA_1969}
and
\citet{Stepanishen_JASA_1971b,Stepanishen_JASA_1971a},
which is an analytic (mesh-free) method dedicated to the evaluation
of transient acoustic fields that are of particular interest
for pulse-echo \gls{us} imaging.

The \gls{sir} approach amounts to evaluating a time-dependent surface integral
derived from the Rayleigh-Sommerfeld equations.
A detailed review of the method can be found in~\cite{Harris_JASA_1981a}.
Once the \gls{sir} is known,
the resulting velocity potential can be computed by the time convolution
of the surface excitation and the \gls{sir}.
Physical quantities such as pressure or particle velocity
can eventually be obtained from the velocity potential.
Much attention and efforts were devoted to the derivation of analytic
\gls{sir} expressions for uniformly excited radiators
of specific shapes and boundary (baffle) conditions,
resulting in complete expressions for
the circular piston~\cite{Harris_JASA_1981a},
the slit~\cite{Lasota_JASA_1984},
the spherically focused radiator~\cite{Arditi_UI_1981},
the triangular piston~\cite{Jensen_JASA_1996},
and the rectangular piston~\cite{SanEmeterio_JASA_1992}.
Semi-analytic expressions were also derived for more complex geometries
such as the cylindrical shell~\cite{Theumann_JASA_1990}
and the toroidal shell~\cite{Baek_JASA_2012}.
Non-uniform excitations were also investigated
in~\cite{Harris_JASA_1981b,Stepanishen_JASA_1981b,Tjotta_JASA_1982,%
Verhoef_JASA_1984},
but analytic expressions were restricted to specific transducer shapes
with specific excitation amplitude distributions.

Analytic expressions only exist for a restricted set of transducer shapes,
excitation distributions, and baffle conditions.
To cope with this limitation,
numerical methods were proposed.
The main principle consists of representing the transducer surface as
a set of characteristic sub-elements for which (simple) analytic expressions
of the \gls{sir} exist.
The \gls{sir} of the transducer is obtained
by summing all the sub-element \glspl{sir}
(superposition principle).
Such a discretization also enables approximating surface apodization and delay
by appropriately weighting and delaying the \gls{sir}
of each sub-element before summation.
A popular characteristic sub-element type is the rectangle,
as it enables reasonably good approximations of arbitrary surfaces,
and more importantly because there exist computationally efficient
far-field approximations for both rigid and soft baffle
conditions~\cite{Jensen_JCA_2001}.
Another relevant approach~\cite{Piwakowski_JASA_1989} consists of
evaluating the surface integral
numerically by
discretizing the surface into ideal points.
The \gls{sir} of the transducer can then be computed by a simple
weighted sum of Dirac delta functions.
This strategy has the advantage of enabling exact representation of surfaces,
at the cost of requiring a much larger number of sub-elements than the strategy
deployed in \gls{fieldii}~\cite{Piwakowski_UFFC_1999}.

A major numerical difficulty of the \gls{sir} method comes from the signal
properties of such responses.
They are typically of very short duration
and are characterized by abrupt slope changes
induced by the edges and vertices of the transducer aperture,
resulting in very high
frequencies~\cite{Lockwood_JASA_1973a,Arditi_UI_1981}.%
\footnote{A typical example of a \gls{sir} and corresponding frequency spectrum
is shown in supplementary material
(\cref{fig:introduction:sir-example}).}
Yet,
because the electromechanical impulse responses of conventional
transducer elements are bandlimited,
the excitation waveform,
and most importantly the field signal of interest (\gls{eg}, pressure field),
are also bandlimited.
Ideally,
one would want to sample the \gls{sir} at the same rate as the excitation
waveform before time convolution of the two signals.
This is obviously not possible in most cases and sampling at the rate of the
\gls{sir} (\gls{ie}, several additional orders of magnitude) is computationally
prohibitive.
Different strategies based on the time integration of the \gls{sir}
were proposed to tackle this computational difficulty,
either numerically using an adaptive sampling
of the \gls{sir}~\cite{Arditi_UI_1981}
or from analytic expressions~\cite{Dhooge_JASA_1997,Jensen_JCA_2001}.
Such strategies enabled for more efficient sampling rates,
but still require oversampling ratios of several factors to prevent from
detrimental aliasing.

In summary,
the most important needs for a computational method relying on the \gls{sir}
are:
an accurate representation of arbitrary transducer shapes;
an efficient way of evaluating the surface integral at each time instant;
and an efficient sampling strategy such that the \gls{sir} can be sampled
at rates similar to those required for the excitation waveform.
To address these needs,
we propose to represent the shape of transducers as \gls{nurbs} surfaces,
to evaluate the surface integral numerically using high-order Gaussian
quadrature rules,
and to express the \gls{sir} in \gls{bspline} bases for efficient sampling
of the time axis.
The use of \gls{nurbs} enables accurate representation of complex surfaces
and exact representations of common transducer elements.
High-order Gaussian quadrature rules
(\gls{eg}, Gauss-Legendre),
require much fewer quadrature points to achieve a desired
accuracy~\cite[Sec.~25.4]{AbramowitzAndStegun_BOOK_1964}.
Finally,
relying on numerical integration allowed us to express the calculation
of the \gls{sir} in high-order \gls{bspline}
bases~\cite{Unser_TSP_1993a,Unser_SPM_1999,Thevenaz_TMI_2000}.
This enabled sampling the \gls{sir} at low sampling rates,
as required by the excitation waveform,
without introducing additional errors on simulated field signals.

%% file: content/approach.tex
\section{Proposed Approach}%
\label{sec:approach}

Let us recall~\cite{Stepanishen_JASA_1971b,Delannoy_JAP_1979b}
that the \gls{sir} of a radiating surface
\( \surface \subset \realnumbers^{3} \)
(assumed to be embedded in an infinite planar baffle)
can be expressed at a field point
\( \fieldpoint \in \volume \subset \realnumbers^{3} \)
and at a time instant \( \timevar \geq 0 \)
as
\begin{align}
    \sirsurface\parens{\fieldpoint, \timevar}
    =
    \int\limits_{\surface}
    \frac{
        \waveformsurfacedistrib\parens{\fieldpointsource}
        \dirac\parens*{
            \timevar
            - \frac{\twonorm{\fieldpoint - \fieldpointsource}}{\soundspeedmean}
        }
    }{
        2 \pi \twonorm{\fieldpoint - \fieldpointsource}
    }
    \sirbcterm\parens{\surfnormalsourceimage, \fieldpoint - \fieldpointsource}
    \ddsurf\parens{\fieldpointsource}
    \label{eq:sir-generic-bcs:bis}
    ,
\end{align}
where \( \soundspeedmean \) is the mean sound speed in the medium
(assumed homogeneous),
\( \surfnormalsourceimage \) is the (inward) surface normal,
\( \waveformsurfacedistrib \) represents the spatial distribution
of surface velocity or pressure over the radiating surface,
and \( \sirbcterm \) is a term depending on the boundary conditions,
expressed as
\begin{align}
    \sirbcterm\parens{\surfnormalsourceimage, \fieldpoint - \fieldpointsource}
    =
    \begin{cases}
        1
        ,
        & \text{rigid baffle}
        ,
        \\
        \cos\parens{\surfnormalsourceimage, \fieldpoint - \fieldpointsource}
        ,
        &
        \text{soft baffle}
        .
    \end{cases}
    \label{eq:sir-generic-bcs-bcterm-def:bis}
\end{align}

Let us consider a generic excitation waveform \( \waveformsurface \)
imposed as boundary condition on the surface \( \surface \)
(\gls{eg}, normal velocity).
The velocity potential can then be expressed as
\(
    \velpot\parens{\fieldpoint, \timevar}
    =
    \waveformsurface\parens{\timevar}
    \conv_{\timevar}
    \sirsurface\parens{\fieldpoint, \timevar}
    ,
\)
from which physical quantities can be obtained.
For instance,
the pressure can be evaluated as
\(
    \pressure\parens{\fieldpoint, \timevar}
    =
    \densitymean
    \parens{\partial / \partial \timevar}
    \velpot\parens{\fieldpoint, \timevar}
\).
Thus,
the acoustic propagation from a radiating surface
can be interpreted as the convolution of the excitation waveform derivative
and the \gls{sir}.
Note that \( \waveformsurface \) is bandlimited
as it includes the electromechanical impulse response of the radiator.

\input{content/approach-quadrature.tex}

\input{content/approach-nurbs.tex}

\input{content/approach-interpolation.tex}

\input{content/approach-implementation.tex}

%% file: content/approach-quadrature.tex
\subsection{Numerical quadrature of the spatial impulse response}%
\label{sec:approach:quadrature}

At each time instant \( \timevar \geq 0 \) and field point
\( \fieldpoint \in \volume \),
the integrand of \cref{eq:sir-generic-bcs:bis}
is a real-valued function that can define as
\(
    \sirsurfaceintegrand\colon
    \surface\subset\realnumbers^{3}\to\realnumbers
\)
for the purpose of derivation.
Let us assume that the radiating surface \( \surface \) is a \emph{smooth}
surface parametrized by a mapping
\(
    \surfacemap \colon
    \surfaceparamspace \subset \realnumbers^{2}
    \to
    \surface
\)
whose \gls{jacobian} determinant
\( \determinant{\surfacemapjac\parens{\surfaceparamvec}} \neq 0 \),
\( \forall \surfaceparamvec \in \surfaceparamspace \),
where the \gls{jacobian} matrix is defined as
\( \surfacemapjac = \parens{\partial / \partial \surfaceparamvec} \surfacemap \).
Thus,
the surface integral of \( \sirsurfaceintegrand \) onto \( \surface \) can be
rewritten as
\begin{align}
    \int\limits_{\surface}
    \sirsurfaceintegrand\parens{\fieldpointsource}
    \ddsurf\parens{\fieldpointsource}
    =
    \int\limits_{\surfaceparamspace}
    \sirsurfaceintegrand\parens{\surfacemap\parens{\surfaceparamvec}}
    \abs{\determinant{\surfacemapjac\parens{\surfaceparamvec}}}
    \ddsurf\parens{\surfaceparamvec}
    .
    \label{eq:surface-integral-mapping}
\end{align}
Assuming that \( \sirsurfaceintegrand \) is a well-behaved function,%
\footnote{This is generally true as there is no realistic interest
for the evaluation of the \gls{sir} on the radiating surface,
where the integrand of \cref{eq:sir-generic-bcs:bis} is obviously singular.}
one can approximate \cref{eq:surface-integral-mapping} by means of Gaussian
quadrature as the finite sum
\begin{align}
    \int\limits_{\surface}
    \sirsurfaceintegrand\parens{\fieldpointsource}
    \ddsurf\parens{\fieldpointsource}
    \approx
    \sum\limits_{\quadpointindex=1}^{\quadpointnb}
    \sirsurfaceintegrand\parens{
        \underbrace{\surfacemap\parens{\surfaceparamvec_{\quadpointindex}}}_{\surfacequadpoint}
    }
    \underbrace{
        \abs{
            \determinant{\surfacemapjac\parens{\surfaceparamvec_{\quadpointindex}}}
        }
    }_{\surfacequadjacdet}
    \surfacequadweight
    ,
    \label{eq:surface-integral-quadrature}
\end{align}
where \( \setroster{\surfacequadcoordvec} \) are the quadrature coordinates
(parametric space),
for which
\( \setroster{\surfacequadweight} \),
\( \setroster{\surfacequadjacdet} \),
and \( \setroster{\surfacequadpoint} \),
are the corresponding quadrature weights,
\gls{jacobian} determinants,
and quadrature points (physical space),
respectively.

Using \cref{eq:surface-integral-quadrature},
we can rewrite \cref{eq:sir-generic-bcs:bis} as
\begin{align}
    \sirsurface\parens{\fieldpoint, \timevar}
    &\approx
    \sum\limits_{\quadpointindex=1}^{\quadpointnb}
    \frac{
        \waveformsurfacedistrib_{\quadpointindex}
        \dirac\parens*{
            \timevar
            - \sirdelay\parens{\fieldpoint, \surfacequadpoint}
        }
    }{
        2 \pi \twonorm{\fieldpoint - \surfacequadpoint}
    }
    \sirbcterm\parens{\surfacequadnormal, \fieldpoint - \surfacequadpoint}
    \surfacequadjacdet
    \surfacequadweight
    \label{eq:sir-generic-bcs-quadrature}
    ,
\end{align}
where
\(
    \sirdelay\parens{\fieldpoint, \surfacequadpoint}
    =
    (1 / \soundspeedmean)
    \twonorm{\fieldpoint - \surfacequadpoint}
\)
and
\(
    \waveformsurfacedistrib_{\quadpointindex}
    =
    \waveformsurfacedistrib\parens{\surfacequadpoint}
\).
The normal vector \( \surfnormal \) at some surface point
\( \fieldpoint = \surfacemap\parens{\surfaceparamvec} \)
can be computed from the surface parametrization (mapping) as
\begin{align}
    \surfnormal
    =
    \frac{
        \partial \surfacemap\parens{\surfaceparamvec}
    }{
        \partial \surfacequadcoordu
    }
    \crossprod
    \frac{
        \partial \surfacemap\parens{\surfaceparamvec}
    }{
        \partial \surfacequadcoordv
    }
    ,
    \label{eq:nurbs-surface-normal}
\end{align}
where \( \parens{\surfacequadcoordu, \surfacequadcoordv} \) are the coordinates
defining the parametric space of the surface mapping.
The \gls{jacobian} determinants,
represented as \( \setroster{\surfacequadjacdet} \) in
\cref{eq:sir-generic-bcs-quadrature},
can then be computed from the normal vectors directly as
\begin{align}
    \abs{\determinant{\surfacemapjac\parens{\surfaceparamvec}}}
    =
    \twonorm*{
        \frac{
            \partial \surfacemap\parens{\surfaceparamvec}
        }{
            \partial \surfacequadcoordu
        }
        \crossprod
        \frac{
            \partial \surfacemap\parens{\surfaceparamvec}
        }{
            \partial \surfacequadcoordv
        }
    }
    .
    \label{eq:nurbs-surface-jacobian-determinant}
\end{align}
Finally,
by grouping all spatially dependent weighting terms in
\cref{eq:sir-generic-bcs-quadrature} under a global weighting term
\( \sirgenericweight \),
we can obtain a compact expression for the evaluation of the \gls{sir}
by means of Gaussian quadrature as
\begin{align}
    \sirsurface\parens{\fieldpoint, \timevar}
    &\approx
    \sum\limits_{\quadpointindex=1}^{\quadpointnb}
    \sirgenericweight\parens{\fieldpoint, \surfacequadpoint}
    \dirac\parens{
        \timevar - \sirdelay\parens{\fieldpoint, \surfacequadpoint}
    }
    \label{eq:sir-generic-bcs-quadrature-compact}
    .
\end{align}

Thus,
provided that such a surface mapping exists and that a suitable Gaussian
quadrature rule is deployed with a sufficient amount of quadrature points
to achieve a desired accuracy,
the calculation of the \gls{sir} can be well approximated
by a sum of shifted-and-weighted Dirac delta functions.

%% file: content/approach-nurbs.tex
\subsection{Non-uniform rational B-spline surface representations}%
\label{sec:approach:nurbs}

\begin{figure*}[htb]
    \centering
    \includegraphics[scale=0.925]{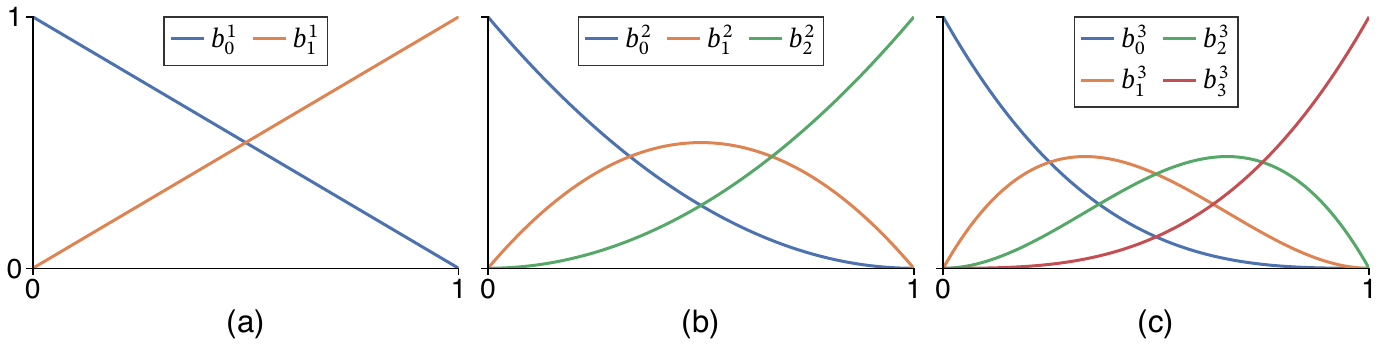}%
    \phantomsubfloatprevspace%
    \phantomsubfloat{\label{fig:approach:nurbs:basis-functions:linear}}%
    \phantomsubfloat{\label{fig:approach:nurbs:basis-functions:quadratic}}%
    \phantomsubfloat{\label{fig:approach:nurbs:basis-functions:cubic}}%
    \caption{%
        (Color online)
        Typical examples of \gls{bspline} basis functions of different degrees
        defined on uniform knot vectors with no interior knots:
        (a)
        linear basis functions defined on the knot vector
        (0, 0, 1, 1);
        (b)
        quadratic basis functions defined on the knot vector
        (0, 0, 0, 1, 1, 1);
        (c)
        cubic basis functions defined on the knot vector
        (0, 0, 0, 0, 1, 1, 1, 1).%
    }%
    \label{fig:approach:nurbs:basis-functions}%
\end{figure*}

We first review some basic principles of \gls{nurbs} (surface) representations.
For a detailed reference,
the reader is referred to the book by
\citet{Piegl_BOOK_1997}.
Let \( \nurbsknotvectoru \) be a nondecreasing sequence of real numbers
representing a nonperiodic, clamped, or open knot vector,
with \( \nurbsknotsnbu + 1 \) knots,
defined as
\begin{align}
    \nurbsknotvectoru
    =
    \tuple{
        \underbrace{\nurbsintervala, \ellipsis, \nurbsintervala}_{\nurbsdegreeu + 1}
        ,
        \nurbscoordu_{\nurbsdegreeu + 1},
        \ellipsis,
        \nurbscoordu_{\nurbsknotsnbu - \nurbsdegreeu - 1},
        \underbrace{\nurbsintervalb, \ellipsis, \nurbsintervalb}_{\nurbsdegreeu + 1}
    }
    ,
\end{align}
\( \nurbsknotsnbu = \nurbsbasisnbu + \nurbsdegreeu + 1 \).
The corresponding (nonnegative) \gls{bspline} basis functions
\( \setroster{\nurbsbasis_{\nurbsindexu}^{\nurbsdegreeu}}_{\nurbsindexu=0}^{\nurbsbasisnbu} \)
of degree \( \nurbsdegreeu \geq 1 \)
are defined recursively as~\cite{Cox_JAM_1972,deBoor_JAT_1972,deBoor_BOOK_1987}
\begin{align}
    \nurbsbasis_{\nurbsindexu}^{0} \parens{\nurbscoordu}
    &=
    \begin{cases}
        1
        ,
        &
        \text{if }
        \nurbscoordu_{\nurbsindexu}
        \leq \nurbscoordu
        < \nurbscoordu_{\nurbsindexu + 1}
        ,
        \\
        0
        ,
        &
        \text{otherwise}
        ,
    \end{cases}
    \\
    \nurbsbasis_{\nurbsindexu}^{\nurbsdegreeu} \parens{\nurbscoordu}
    &=
    \frac{
        \nurbscoordu - \nurbscoordu_{\nurbsindexu}
    }{
        \nurbscoordu_{\nurbsindexu + \nurbsdegreeu} - \nurbscoordu_{\nurbsindexu}
    }
    \nurbsbasis_{\nurbsindexu}^{\nurbsdegreeu - 1} \parens{\nurbscoordu}
    +
    \frac{
        \nurbscoordu_{\nurbsindexu + \nurbsdegreeu + 1} - \nurbscoordu
    }{
        \nurbscoordu_{\nurbsindexu + \nurbsdegreeu + 1} - \nurbscoordu_{\nurbsindexu + 1}
    }
    \nurbsbasis_{\nurbsindexu + 1}^{\nurbsdegreeu - 1} \parens{\nurbscoordu}
    .
    \label{eq:bspline-basis-recursive}
\end{align}
Note that we restrict ourselves to nonperiodic knot vectors onto which
\gls{bspline} basis functions are interpolating at the endpoints of
such knot vectors but are (in general) non-interpolating at interior knots.
Also,
definitions are not strictly limited to the
\( \bracks{\nurbsintervala, \nurbsintervalb} \) interval.
Yet, it is so common to the \gls{nurbs} community that it is also adopted
in the present work.
Typical examples of \gls{bspline} basis functions on knot vectors
of the form
\(
    \nurbsknotvectoru
    =
    \tuple{
        \nurbsintervala, \ellipsis, \nurbsintervala,
        \nurbsintervalb, \ellipsis, \nurbsintervalb
    }
\)
are shown in \cref{fig:approach:nurbs:basis-functions}
for different degrees.
This type of knot vectors results in \gls{bspline} basis functions
that are \gls{bernstein} polynomials,
because \gls{bspline} representations are a generalization
of \gls{bezier} representations~\cite[Sec.~2.2]{Piegl_BOOK_1997}.

A \gls{nurbs} surface \( \surface \)
of degree \( \parens{\nurbsdegreeu, \nurbsdegreev} \)
in directions \( \parens{\nurbscoordu, \nurbscoordv} \)
can be represented by a bivariate vector-valued piecewise rational function
(mapping)
\(
    \nurbssurfacemap\colon
    \bracks{\nurbsintervala, \nurbsintervalb}^{2}
    \to
    \nurbssurface\subset\realnumbers^{3}
\)
defined as~\cite[Sec.~4.4]{Piegl_BOOK_1997}
\begin{align}
    \nurbssurfacemap\parens{\nurbscoordu, \nurbscoordv}
    =
    \dfrac{
        \displaystyle  %
        \sum\limits_{\nurbsindexu = 0}^{\nurbsbasisnbu}
        \sum\limits_{\nurbsindexv = 0}^{\nurbsbasisnbv}
        \nurbsbasis_{\nurbsindexu}^{\nurbsdegreeu}\parens{\nurbscoordu}
        \nurbsbasis_{\nurbsindexv}^{\nurbsdegreev}\parens{\nurbscoordv}
        \nurbsweight_{\nurbsindexu, \nurbsindexv}
        \nurbsctrlpoint_{\nurbsindexu, \nurbsindexv}
    }{
        \displaystyle  %
        \sum\limits_{\nurbsindexu = 0}^{\nurbsbasisnbu}
        \sum\limits_{\nurbsindexv = 0}^{\nurbsbasisnbv}
        \nurbsbasis_{\nurbsindexu}^{\nurbsdegreeu}\parens{\nurbscoordu}
        \nurbsbasis_{\nurbsindexv}^{\nurbsdegreev}\parens{\nurbscoordv}
        \nurbsweight_{\nurbsindexu, \nurbsindexv}
    }
    ,
    \label{eq:nurbs-surface-mapping}
\end{align}
where
\( \setroster{\nurbsctrlpoint_{\nurbsindexu, \nurbsindexv}} \)
are the control points
(forming a bidirectional net),
\( \setroster{\nurbsweight_{\nurbsindexu, \nurbsindexv}} \)
are the corresponding weights,
\( \setroster{\nurbsbasis_{\nurbsindexu}^{\nurbsdegreeu}} \)
and
\( \setroster{\nurbsbasis_{\nurbsindexv}^{\nurbsdegreev}} \)
are the (nonrational) \gls{bspline} basis functions
of degrees \( \nurbsdegreeu \) and \( \nurbsdegreev \)
that are defined on the nondecreasing and nonperiodic
(\gls{ie}, clamped)
knot vectors
\begin{align}
    \nurbsknotvectoru
    &=
    \tuple{
        \underbrace{\nurbsintervala, \ellipsis, \nurbsintervala}_{\nurbsdegreeu + 1}
        ,
        \nurbscoordu_{\nurbsdegreeu + 1},
        \ellipsis,
        \nurbscoordu_{\nurbsknotsnbu - \nurbsdegreeu - 1},
        \underbrace{\nurbsintervalb, \ellipsis, \nurbsintervalb}_{\nurbsdegreeu + 1}
    }
    ,
    \\
    \nurbsknotvectorv
    &=
    \tuple{
        \underbrace{\nurbsintervala, \ellipsis, \nurbsintervala}_{\nurbsdegreev + 1}
        ,
        \nurbscoordv_{\nurbsdegreev + 1},
        \ellipsis,
        \nurbscoordv_{\nurbsknotsnbv - \nurbsdegreev - 1},
        \underbrace{\nurbsintervalb, \ellipsis, \nurbsintervalb}_{\nurbsdegreev + 1}
    }
    ,
\end{align}
respectively,
where
\( \nurbsknotsnbu = \nurbsbasisnbu + \nurbsdegreeu + 1 \)
and
\( \nurbsknotsnbv = \nurbsbasisnbv + \nurbsdegreev + 1 \).
Expressions for the derivatives of \gls{nurbs} surfaces
also exist and can be found in~\cite[Sec.~4.5]{Piegl_BOOK_1997}.
They are essential tools for the calculations of the surface normal vectors
and the \gls{jacobian} determinants defined in
\cref{eq:nurbs-surface-normal} and \cref{eq:nurbs-surface-jacobian-determinant},
respectively.
Note that if both \( \nurbsknotvectoru \) and \( \nurbsknotvectorv \)
are defined with no interior points,
such as in the examples shown
in~\cref{fig:approach:nurbs:basis-functions},
the \gls{nurbs} surface is a (smooth) rational \gls{bezier} surface.

An important property of \gls{nurbs} representations is that any \gls{nurbs}
surface can be decomposed into a union of (smooth) rational \gls{bezier} surfaces
(or patches),
that share at most a common edge or a common vertex.
Such a decomposition can be performed by a simple procedure known as
knot refinement~\cite[Sec.~5.3]{Piegl_BOOK_1997}.
Each (decomposed) rational \gls{bezier} patch is actually also
a \gls{nurbs} surface such that it can be readily
represented by \cref{eq:nurbs-surface-mapping}.
Rational \gls{bezier} patches are of particular interest as they are
\emph{smooth} surfaces as opposed to \enquote{general} \gls{nurbs} surfaces
that may contain breakpoints.
This is an essential property to allow the use of Gaussian quadrature rules
onto such surfaces
(\cref{sec:approach:quadrature}).

Also,
most transducer elements composing conventional \gls{us} transducer arrays
are (at most) quadric surfaces
and \gls{nurbs} surfaces can represent quadric surfaces
exactly
(suppl.\@ mat., \cref{sec:sup:approach:nurbs}).

%% file: content/approach-interpolation.tex
\subsection{Spatial impulse response in B-spline bases}%
\label{sec:approach:interpolation}

Thanks to the \gls{nurbs} representation of the radiating surface,
we can evaluate the global weigthing term
\( \sirgenericweight \) from the sum of (shifted-and-weighted)
Dirac delta functions to approximate the \gls{sir}
[\cref{eq:sir-generic-bcs-quadrature-compact}].
Before diving into the well-known difficulty of sampling Dirac delta functions,
it is important to keep in mind that the \gls{sir} is mainly a
physical concept as such a quantity cannot be measured in physical conditions
due to the electromechanical impulse response of transducer elements
that are bandlimited.
Thus,
signals of interest that will be eventually measured are also bandlimited.
Such a signal can be expressed generically as
\begin{align}
    \sirgenericsignal\parens{\fieldpoint, \timevar}
    =
    \sirgenericexcitation\parens{\timevar}
    \conv_{\timevar}
    \sirsurface\parens{\fieldpoint, \timevar}
    \label{eq:sir-generic-bandpass-signal}
    ,
\end{align}
where \( \sirgenericexcitation \) represents some bandpass excitation waveform.
Using \cref{eq:sir-generic-bcs-quadrature-compact} we obtain the corresponding
approximation of such a signal as
\begin{align}
    \sirgenericsignal\parens{\fieldpoint, \timevar}
    \approx
    \sirgenericexcitation\parens{\timevar}
    \conv_{\timevar}
    \sum\limits_{\quadpointindex=1}^{\quadpointnb}
    \sirgenericweight\parens{\fieldpoint, \surfacequadpoint}
    \dirac\parens{
        \timevar - \sirdelay\parens{\fieldpoint, \surfacequadpoint}
    }
    \label{eq:sir-generic-bandpass-signal-quadrature}
    .
\end{align}
\Cref{eq:sir-generic-bandpass-signal-quadrature} tells us that the signal
\( \sirgenericsignal \) is composed of shifted-and-weighted replicas
of the excitation waveform \( \sirgenericexcitation \).
The excitation waveform \( \sirgenericexcitation \) is typically measured
in physical conditions or approximated using a model pulse.
Thus,
it can be expressed as a sampled signal of
uniformly spaced samples
\(
    \setroster{
        \sirgenericexcitation\parens{\sirexcitationtimeindex\sirdeltat}
    }
\),
\(
    \sirexcitationtimeindex
    \in
    \bracks{0, \dots, \sirexcitationtimenb - 1}
\),
with a sampling interval \( \sirdeltat \).
\Cref{eq:sir-generic-bandpass-signal-quadrature}
can then be rewritten as the discrete convolution
\begin{align}
    \sirgenericsignal\parens{\fieldpoint, \timevar}
    &\approx
    \sum\limits_{\sirexcitationtimeindex \in \integernumbers}
    \sirgenericexcitation\parens{\sirexcitationtimeindex\sirdeltat}
    \sum\limits_{\quadpointindex=1}^{\quadpointnb}
    \frac{
        \sirgenericweight\parens{\fieldpoint, \surfacequadpoint}
    }{
        \sirdeltat
    }
    \interpbasis\parens*{
        \frac{
            \timevar - \sirdelay\parens{\fieldpoint, \surfacequadpoint}
        }{
            \sirdeltat
        }
        - \sirexcitationtimeindex
    }
    \label{eq:sir-generic-bandpass-signal-quadrature-discrete-fast}
    ,
\end{align}
where \( \interpbasis \) is some (interpolating) basis function
(sometimes referred to as sampling kernel),
and where it is assumed that boundary conditions are handled properly
(via signal extension of \( \sirgenericexcitation \)).
Since \( \sirgenericexcitation \) is a bandpass signal,
a natural and error-free choice for \( \interpbasis \) would be the (normalized)
\gls{sinc} function,
\( \sinc\parens{x} \coloneqq \sin\parens{\pi x} / (\pi x) \).
But because the \gls{sinc} function is of infinite support,
it needs to be truncated
by multiplying it by some finite-support window.
Generally,
truncated \gls{sinc} functions require rather large supports to achieve
high accuracy~\cite{Thevenaz_TMI_2000}.
Other options are standard (finite-support) interpolating basis functions
such as those deployed for nearest-neighbor and linear interpolation,
which are computationally efficient but of low accuracy.

Instead,
as the discrete convolution of
\cref{eq:sir-generic-bandpass-signal-quadrature-discrete-fast}
consists of shifting and interpolating \( \sirgenericexcitation \),
we propose to rely on the concept of generalized interpolation~\cite{Thevenaz_TMI_2000}.
By doing so,
we can rewrite
\cref{eq:sir-generic-bandpass-signal-quadrature-discrete-fast}
as
\begin{align}
    \sirgenericsignal\parens{\fieldpoint, \timevar}
    &\approx
    \sum\limits_{\sirexcitationtimeindex \in \integernumbers}
    \sirgenericexcitationcoeffs\parens{\sirexcitationtimeindex\sirdeltat}
    \sum\limits_{\quadpointindex=1}^{\quadpointnb}
    \frac{
        \sirgenericweight\parens{\fieldpoint, \surfacequadpoint}
    }{
        \sirdeltat
    }
    \interpbasisgeneral\parens*{
        \frac{
            \timevar - \sirdelay\parens{\fieldpoint, \surfacequadpoint}
        }{
            \sirdeltat
        }
        - \sirexcitationtimeindex
    }
    \label{eq:sir-generic-bandpass-signal-quadrature-discrete-geninterp-fast}
    ,
\end{align}
where \( \interpbasisgeneral \) is some basis function
and
\(
    \setroster{
        \sirgenericexcitationcoeffs\parens{\sirexcitationtimeindex\sirdeltat}
    }
\)
are the corresponding basis coefficients defined such that
\begin{align}
    \sirgenericexcitation\bracks{\sirexcitationtimeindex_{0}}
    =
    \sum\limits_{\sirexcitationtimeindex \in \integernumbers}
    \sirgenericexcitationcoeffs\bracks{\sirexcitationtimeindex}
    \interpbasisgeneral\parens{
        \sirexcitationtimeindex_{0} - \sirexcitationtimeindex
    }
    .
\end{align}
As a result,
one can rely on potentially non-interpolating basis functions,
possessing better properties than interpolating ones,
with an additional operation that consists of finding the basis coefficients.
A particularly elegant and efficient way of finding these coefficients is
yet another (pre-filtering) convolution operation with the convolution-inverse%
\footnote{Not to be confused with the inverse function of
\( \interpbasisgeneral \).}
\( \interpbasisconvinv \) as
\begin{align}
    \sirgenericexcitationcoeffs\bracks{\sirexcitationtimeindex_{0}}
    =
    \sum\limits_{\sirexcitationtimeindex \in \integernumbers}
    \interpbasisconvinv\parens{\sirexcitationtimeindex_{0}}
    \sirgenericexcitation\bracks{
        \sirexcitationtimeindex_{0} - \sirexcitationtimeindex
    }
    ,
    \label{eq:interp-pre-filtering}
\end{align}
provided that \( \interpbasisconvinv \) exists.
In the case of interpolating basis functions,
the sequence of coefficients
\(
    \setroster{
        \sirgenericexcitationcoeffs\parens{\sirexcitationtimeindex\sirdeltat}
    }
\)
is equal to the sequence of samples
\(
    \setroster{
        \sirgenericexcitation\parens{\sirexcitationtimeindex\sirdeltat}
    }
\).

We will restrict ourselves to the family of \gls{bspline} basis functions
as they benefit from maximal approximation orders for given
supports~\cite{Unser_TSP_1993a,Unser_SPM_1999,Thevenaz_TMI_2000}.
For the purpose of interpolation,
we can express a \gls{bspline} basis function of degree
\( \interpbsplinedegree \) as~\cite{Schoenberg_QAM_1946a,Schoenberg_QAM_1946b}
\begin{align}
    \interpbspline^{\interpbsplinedegree}\parens{x}
    =
    \sum\limits_{k=0}^{\interpbsplinedegree+1}
    \frac{
        \parens{-1}^{k}
        \parens{\interpbsplinedegree + 1}
    }{
        \parens{\interpbsplinedegree+1-k}!\,k!
    }
    \parens*{
        \frac{\interpbsplinedegree + 1}{2} + x - k
    }_{+}^{\interpbsplinedegree}
    ,
\end{align}
\( \forall x \in \realnumbers \),
\( \forall \interpbsplinedegree \in \realnumbers \),
where the one-sided power function \( \parens{x}_{+}^{\interpbsplinedegree} \)
is defined as~\cite{Thevenaz_TMI_2000}
\begin{align}
    \parens{x}_{+}^{\interpbsplinedegree}
    =
    \begin{cases}
        0
        & \interpbsplinedegree = 0 \text{ and } x < 0
        ,
        \\
        1/2
        & \interpbsplinedegree = 0 \text{ and } x = 0
        ,
        \\
        1
        & \interpbsplinedegree = 0 \text{ and } x > 0
        ,
        \\
        \parens{x}_{+}^{0} x^{\interpbsplinedegree}
        & \interpbsplinedegree >0
        .
    \end{cases}
\end{align}
We are typically interested in \gls{bspline} basis functions
\( \interpbspline^{\interpbsplinedegree} \) of degree greater than one,
since \( \interpbspline^{0} \) is almost identical to the nearest-neighbor
basis function and \( \interpbspline^{1} \) corresponds to the basis function
of linear interpolation.
Also,
\gls{bspline} basis functions of even degrees are usually not
computationally interesting because they have the same support
as \gls{bspline} basis functions of the following (odd) degrees.
Since \gls{bspline} basis functions are symmetric,
the pre-filtering operation defined in \cref{eq:interp-pre-filtering}
can be efficiently performed by a series of
\( \interpbsplinepolepairnb = \floor{n/2} \) consecutive
causal and anti-causal \gls{iir} filters
with poles
\( \setroster{\interpbsplinepole_{\interpbsplinepoleindex}} \)
and
\( \setroster{\interpbsplinepole\inv_{\interpbsplinepoleindex}} \),
respectively~\cite{Unser_PAMI_1991,Unser_TSP_1993a,Unser_TSP_1993b}.
The \( \interpbsplinepolepairnb \) pairs of poles can be derived
from the \gls{z-transform} of the convolution-inverse
\( \parens{\interpbspline^{\interpbsplinedegree}}\inv \),
many of which are tabulated
in~\cite{Unser_PAMI_1991,Unser_TSP_1993a,Unser_TSP_1993b}.

So far we have assumed that the samples
\(
    \setroster{
        \sirgenericexcitation\parens{\sirexcitationtimeindex\sirdeltat}
    }
\)
and corresponding coefficients
\(
    \setroster{
        \sirgenericexcitationcoeffs\parens{\sirexcitationtimeindex\sirdeltat}
    }
\)
were extended properly
\(
    \forall
    \sirexcitationtimeindex
    \in
    \integernumbers
    \setminus
    \bracks{0, \dots, \sirexcitationtimenb - 1}
\)
for the purpose of convolution operations.
Because the corresponding excitation waveform is bandlimited and can
be modeled as a windowed sinusoidal \gls{rf} signal of
finite support,%
\footnote{Acoustic pulses are theoretically infinite responses with
(rapidly) decaying trailing oscillations.
Thus,
they may be truncated when their trailing pulse envelope falls below
some level to achieve a desired accuracy.}
we consider zero boundary conditions.
Also,
since we want to work with finite support sequences to perform
the discrete convolution
in~\cref{eq:sir-generic-bandpass-signal-quadrature-discrete-geninterp-fast},
we impose zero boundary conditions to the coefficients directly,
namely
\(
    \sirgenericexcitationcoeffs\parens{\sirexcitationtimeindex\sirdeltat}
    =
    0
\),
\(
    \forall
    \sirexcitationtimeindex
    \in
    \integernumbers
    \setminus
    \bracks{0, \dots, \sirexcitationtimenb - 1}
\).

Once these coefficients are computed,
they can be used readily for all quadrature points to evaluate
the signal of interest \( \sirgenericsignal \) at any field point
\( \fieldpoint \) using
\cref{eq:sir-generic-bandpass-signal-quadrature-discrete-geninterp-fast},
namely as
a discrete (full) convolution of the (pre-filtered) coefficients
and the summation of \( \quadpointnb \) shifted-and-weighted
basis functions
\gls{wrt}
the propagation times
\( \setroster{\sirdelay\parens{\fieldpoint, \surfacequadpoint}} \)
and weights
\( \setroster{\sirgenericweight\parens{\fieldpoint, \surfacequadpoint}} \),
respectively.
An illustration of the convolution involved in the evaluation of
\cref{eq:sir-generic-bandpass-signal-quadrature-discrete-geninterp-fast}
is described in supplementary material
(\cref{sec:sup:approach:implementation}).
For the purpose of interpretation (and compactness),
let us define the basis \gls{sir} as
\begin{align}
    \sirsurfacegeneral\parens{\fieldpoint, \timevar}
    =
    \sum\limits_{\quadpointindex=1}^{\quadpointnb}
    \frac{
        \sirgenericweight\parens{\fieldpoint, \surfacequadpoint}
    }{
        \sirdeltat
    }
    \interpbasisgeneral\parens*{
        \frac{
            \timevar - \sirdelay\parens{\fieldpoint, \surfacequadpoint}
        }{
            \sirdeltat
        }
        - \sirexcitationtimeindex
    }
    \label{eq:basis-sir-general}
    ,
\end{align}
such that
\cref{eq:sir-generic-bandpass-signal-quadrature-discrete-geninterp-fast}
can be interpreted as the convolution of the (pre-filtered) coefficients
and the basis \gls{sir},
namely
\begin{align}
    \sirgenericsignal\parens{\fieldpoint, \timevar}
    &\approx
    \sum\limits_{\sirexcitationtimeindex \in \integernumbers}
    \sirgenericexcitationcoeffs\parens{\sirexcitationtimeindex\sirdeltat}
    \sirsurfacegeneral\parens*{
        \fieldpoint,
        \timevar - \sirexcitationtimeindex\sirdeltat
    }
    \label{eq:sir-generic-bandpass-signal-quadrature-discrete-geninterp-fast-compact}
    .
\end{align}

Even if \( \sirsurface \) is not a quantity that can be measured,
and as such is less important to approximate accurately than
\( \sirgenericsignal \),
it remains interesting to be able to evaluate it.%
\footnote{Note that it is also possible to obtain the field response to
a continuous-wave excitation directly from the \gls{sir}
by evaluating its \gls{fourier} transform.}
Since \( \sirsurfacegeneral \) is a linear combination
of \( \quadpointnb \) basis functions
[\cref{eq:basis-sir-general}],
it is possible to obtain an approximation of \( \sirsurface \)
from \( \sirsurfacegeneral \)
with the convolution-inverse \( \interpbasisconvinv \)
as
\begin{align}
    \sirsurface\parens{\fieldpoint, \timevar}
    &\approx
    \sum\limits_{\sirexcitationtimeindex \in \integernumbers}
    \interpbasisconvinv\parens{\sirexcitationtimeindex}
    \sirsurfacegeneral\parens*{
        \fieldpoint,
        \timevar - \sirexcitationtimeindex\sirdeltat
    }
    \label{eq:sir-from-basis-sir-conv-inverse}
    .
\end{align}
This is similar to expressing a cardinal spline basis function from
its corresponding (non-interpolating) \gls{bspline} basis function.
Thus,
the approximation of the \gls{sir} obtained from
\cref{eq:sir-from-basis-sir-conv-inverse}
when considering non-interpolating basis functions will be of infinite
support with rapidly decaying oscillations~\cite{Thevenaz_TMI_2000}.

%% file: content/approach-implementation.tex
\subsection{Complete process and implementation details}%
\label{sec:approach:implementation}

\begin{figure*}[htb]
    \centering
    \includegraphics[scale=0.925]{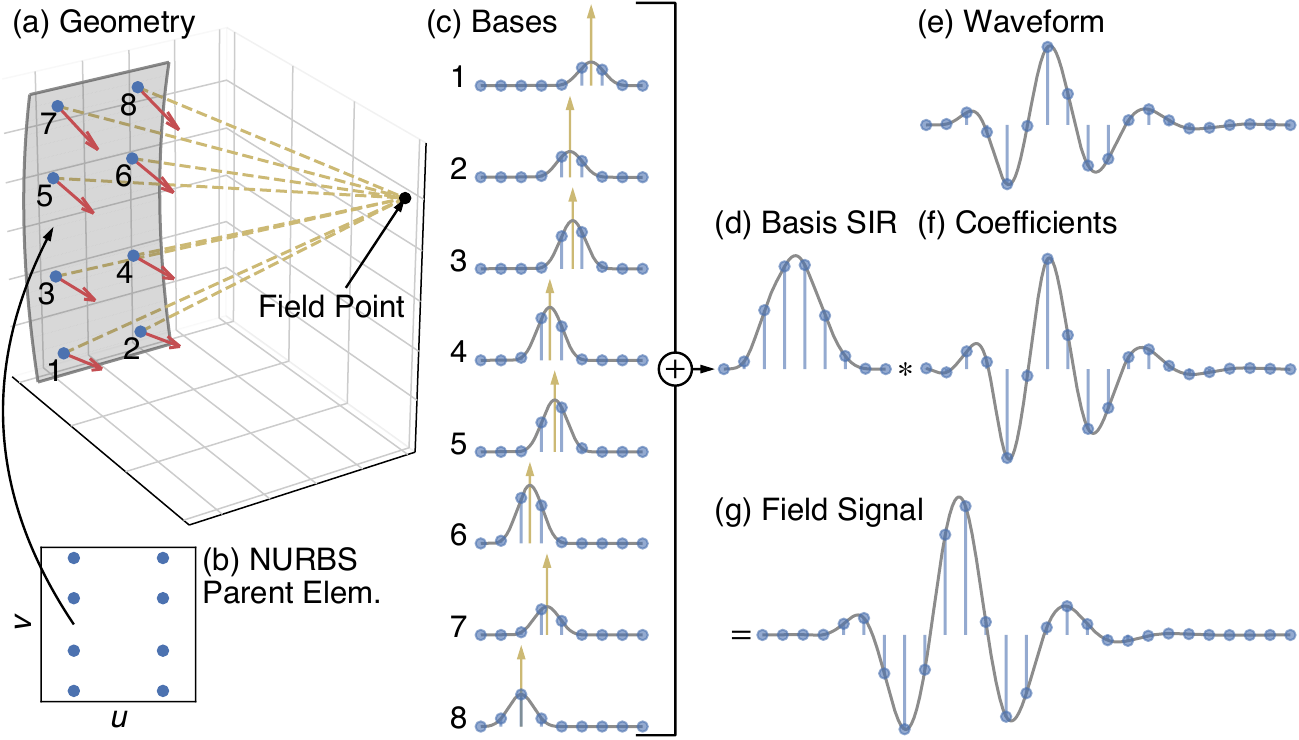}%
    \phantomsubfloatprevspace%
    \phantomsubfloat{\label{fig:approach:complete-strategy:geom}}%
    \phantomsubfloat{\label{fig:approach:complete-strategy:nurbs}}%
    \phantomsubfloat{\label{fig:approach:complete-strategy:bases}}%
    \phantomsubfloat{\label{fig:approach:complete-strategy:basis-sir}}%
    \phantomsubfloat{\label{fig:approach:complete-strategy:waveform}}%
    \phantomsubfloat{\label{fig:approach:complete-strategy:coeffs}}%
    \phantomsubfloat{\label{fig:approach:complete-strategy:signal}}%
    \caption{%
        (Color online)
        Illustration of all steps composing the proposed approach
        for the approximation of field signals using the \glsxtrfull{sir} of
        radiating surfaces.
        (Please refer to the associated text for a summary of each step
        involved.)%
    }%
    \label{fig:approach:complete-strategy}
\end{figure*}

\Cref{fig:approach:complete-strategy} summarizes
the complete process of the proposed approach developed in
\cref{sec:approach:quadrature,sec:approach:nurbs,%
sec:approach:interpolation}
for the approximation of (bandpass) field signals using the \gls{sir} of
radiating surfaces.
The cylindrical shape of the transducer element
[\cref{fig:approach:complete-strategy:geom}]
is represented (exactly) as a \gls{nurbs} surface
using \cref{eq:nurbs-surface-mapping}.
As this \gls{nurbs} surface is also a \emph{smooth} \gls{bezier} surface,
a \( \parens{2 \times 4} \) Gauss-Legendre quadrature rule
can be defined in the parametric space directly
[\cref{fig:approach:complete-strategy:nurbs}],
The corresponding eight quadrature points can be mapped into the physical
space using the \gls{nurbs} representation.
The normal vectors and \gls{jacobian} determinants
can be computed using \cref{eq:nurbs-surface-normal}
and \cref{eq:nurbs-surface-jacobian-determinant},
respectively.
From the distances between the quadrature points and the field point,
some basis function (\gls{eg}, cubic \gls{bspline}) can be evaluated
at the corresponding time instants on a sampled time axis of minimum support
[\cref{fig:approach:complete-strategy:bases}],
and weighted accordingly.
Their summation
[\cref{eq:basis-sir-general}]
results in the basis \gls{sir}
[\cref{fig:approach:complete-strategy:basis-sir}].
The (pre-filtered) coefficients
[\cref{fig:approach:complete-strategy:coeffs}]
are obtained by a series of causal and anti-causal \gls{iir} filters
applied to the excitation waveform.
Finally,
the field signal
[\cref{fig:approach:complete-strategy:signal}]
is obtained by the convolution of
the (pre-filtered) coefficients [\cref{fig:approach:complete-strategy:coeffs}]
and
the basis \gls{sir} [\cref{fig:approach:complete-strategy:basis-sir}].
This process must be repeated for all field points of interest,
except for the coefficients that only need to be computed once.

The extension of the proposed approach to transducer arrays and pulse-echo
acquisitions is straightforward and can be efficiently implemented
(suppl.\@ mat., \cref{sec:approach:arrays}).

%% file: content/validation.tex
\section{Experiments and Results}%
\label{sec:validation}

To validate the proposed approach,
we performed two numerical experiments
(\cref{sec:validation:interpolation,%
sec:validation:analytic}).
The goal of the first one is to validate the core of the proposed approach,
as described in
\cref{eq:sir-generic-bandpass-signal-quadrature-discrete-geninterp-fast-compact},
namely a convolution of (pre-filtered) coefficients and a signal
(representing the basis \gls{sir}) composed of shifted-and-weighted
basis functions.
The second experiment consists of evaluating the accuracy
of the complete proposed approach on field signals radiated by transducer elements
with specific shapes allowing analytic expressions for the \gls{sir}.

For the two experiments,
the pulse model considered
was an analytic expression of the time derivative
of a log-normal-modulated sinusoidal \gls{rf} pulse
(suppl.\@ mat., \cref{sec:sup:validation}).
For the log-normal distribution
\( \lognormaldist\parens{\mu, \sigma} \),
we used the parameters
\(
    \mu
    \approx
    \num[round-mode=places,round-precision=2]{-14.802665843192596}
\)
and
\(
    \sigma
    \approx
    \num[round-mode=places,round-precision=2]{0.25551258495020857}
\).
The frequency of the sinusoidal was set to \SI{4.75}{\mega\hertz},
resulting in a waveform centered at \( \SI{\sim 5.3}{\mega\hertz} \)
with a bandwidth of \( \SI{\sim 71}{\percent} \) at \SI{-6}{\decibel}
[\cref{fig:validation:waveform:spectrum}].
The waveform was truncated at a trailing pulse envelope level
of \SI{-320}{\decibel} (double-precision floating-point format),
resulting in a duration of \( \SI{\sim 3.15}{\micro\second} \).
The \gls{fwhm} of the resulting waveform
corresponds to
\( \SI{\sim 0.23}{\micro\second} \).

In addition to the \gls{bspline} basis functions,
we also considered \gls{o-moms} functions~\cite{Blu_TIP_2001},
that are derived from \gls{bspline} basis functions.
In general,
we compared the classical interpolating basis functions
for nearest-neighbor,
linear,
and
quadratic%
\footnote{By quadratic,
it is implied that the optimal \( a = -1/2 \) parameter was used
for the \gls{keys} basis function,
resulting in a quadratic interpolation method of order three.}
\gls{keys}~\cite{Keys_TASSP_1981},
as well as \gls{bspline} and \gls{o-moms} basis functions of different degrees.
For the purpose of quantifying the accuracy in the following experiments,
we relied on the relative two-norm error
defined between an estimated (\gls{rf} signal) vector and its true counterpart
\(
    \estimate{\sirgenericsignalvec}, \sirgenericsignalvec
    \in \realnumbers^{n}
\)
as
\begin{align}
    \varepsilon
    =
    \twonorm{\sirgenericsignalvec - \estimate{\sirgenericsignalvec}}
    /
    \twonorm{\sirgenericsignalvec}
    \label{eq:sir-valid-exp-rel-error}
    .
\end{align}

\input{content/validation-interpolation.tex}

\input{content/validation-analytic.tex}

%% file: content/validation-interpolation.tex
\subsection{Convergence order of various basis functions}%
\label{sec:validation:interpolation}

We considered a reference signal consisting of the time convolution between
the analytic waveform \( \sirexpwaveform \) considered
[\cref{fig:validation:waveform:shape}]
and a stream of \( n \) Dirac delta functions
at random-uniform times
\(
    \setroster{
        \sirexpdiracdelay_{\sirexpdiracindex}
    }
\)
and of random-normal amplitudes
\(
    \setroster{
        \sirexpdiracweight_{\sirexpdiracindex}
    }
\).
The corresponding reference signal can be expressed as
\begin{align}
    \sirexpsignal\parens{\timevar}
    =
    \sirexpwaveform\parens{\timevar}
    \conv_{\timevar}
    \sum\limits_{\sirexpdiracindex=1}^{\sirexpdiracnb}
    \sirexpdiracweight_{\sirexpdiracindex}
    \dirac\parens{
        \timevar - \sirexpdiracdelay_{\sirexpdiracindex}
    }
    \label{eq:sir-exp-convergence-ref-signal}
    ,
\end{align}
which can be seen as a generic expression for
\cref{eq:sir-generic-bandpass-signal-quadrature},
namely signals we want to approximate using the proposed approach
\cref{eq:sir-generic-bandpass-signal-quadrature-discrete-geninterp-fast-compact}.
Note that \cref{eq:sir-exp-convergence-ref-signal}
can be evaluated exactly using the analytic expression for the
excitation waveform.

We considered a signal duration corresponding to \num{500} times
the time resolution cell of the waveform considered
and populated it with an average of \num{100} random Diracs
per time resolution cell.
We compared three interpolating basis functions,
namely nearest-neighbor (degree zero), linear (degree one),
and quadratic \gls{keys} (degree two),
four (non-interpolating) \gls{bspline} basis functions of degrees
\numlist[]{2;3;4;5},
and the (non-interpolating) \gls{o-moms} basis function of degree three.
To validate that the theoretical approximation order
of these basis functions (\gls{ie}, degree plus one)
is achieved by the proposed method,
we compared each approximated signal with the analytic one
using \cref{eq:sir-valid-exp-rel-error}
at \num{15} log-linearly spaced sampling rates
ranging from \SI{20}{\mega\hertz} to \SI{1}{\giga\hertz}.

\begin{figure}[htb]
    \centering
    \includegraphics[scale=0.925]{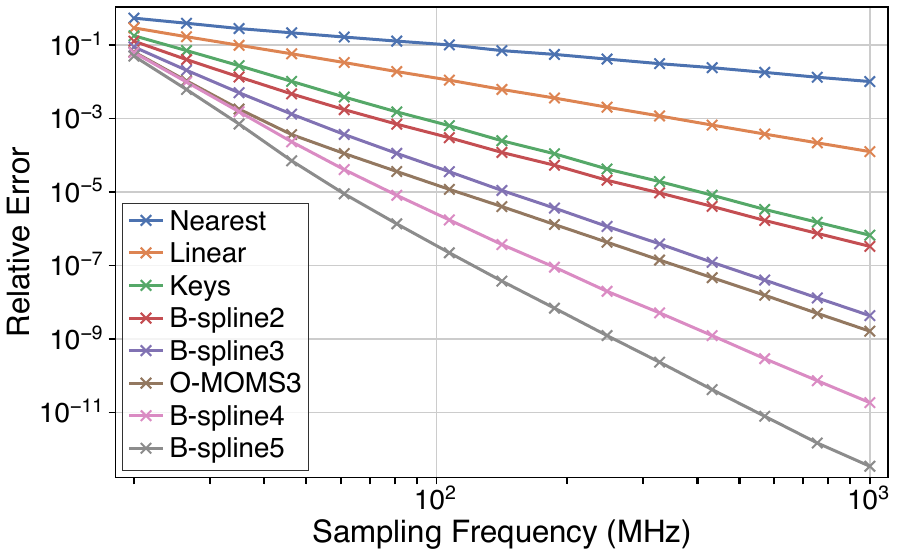}%
    \caption{%
        (Color online)
        Results of the numerical experiment performed to validate
        the theoretical orders of convergence for different basis functions.
        Different types of basis functions were considered:
        three interpolating basis functions,
        namely nearest-neighbor (degree zero),
        linear (degree one),
        and quadratic \gls{keys} (degree two);
        four non-interpolating \gls{bspline} basis functions of degree
        \numlist[]{2;3;4;5};
        and
        the non-interpolating \glsxtrfull{o-moms} basis function
        of degree three.
        The theoretical convergence order of a basis function
        is equal to its degree plus one.%
    }%
    \label{fig:validation:interpolation:convergence-order}
\end{figure}

The resulting convergence curves are depicted in
\cref{fig:validation:interpolation:convergence-order}.
One can see that the theoretical approximation orders are accurately
validated for all basis functions
and that the use of high-order non-interpolating \gls{bspline}
basis functions provides a major advantage.
It can also be mentioned that the \gls{o-moms} of degree three seems to be
\enquote{over-performing} at low frequencies,
namely from \SIrange[]{20}{50}{\mega\hertz}.
This low-frequency range represents a \enquote{rough} Nyquist-rate range
for the excitation waveform considered
[\cref{fig:validation:waveform:spectrum}],
especially below \SI{30}{\mega\hertz}.

%% file: content/validation-analytic.tex
\subsection{Validation against analytic solutions}%
\label{sec:validation:analytic}

We considered two transducer-element shapes and corresponding baffle conditions
for which analytic expressions are available for evaluating the \gls{sir}
at any field point:
the spherically focused element (\gls{ie}, spherical cap)
with a rigid baffle boundary~\cite{Arditi_UI_1981}
and
the rectangular plane
with a soft baffle condition~\cite{SanEmeterio_JASA_1992}.
In all cases,
we considered the same excitation waveform
[in \cref{fig:validation:waveform}],
characterized by a (center) wavelength
\( \wavelength \approx \SI{291}{\micro\meter} \)
for \( \soundspeedmean = \SI{1540}{\meter\per\second} \).
We evaluated the field signal radiated by the transducer element at three
characteristic field points (A, B, and C)
positioned relatively to the transducer element center.
The three field points were positioned such that their projections onto the
surface lie on an axis of symmetry, an edge, and outside the surface projection.

We restricted the comparison to the following basis functions:
nearest-neighbor, linear, quadratic \gls{keys},
\gls{bspline} of degrees three and five,
and \gls{o-moms} of degree three.
Two sampling rates of \SI{30}{\mega\hertz} and \SI{80}{\mega\hertz}
were considered,
namely a rather low sampling rate for the excitation waveform considered,
and a rather high one
[\cref{fig:validation:waveform:spectrum}].
Both shapes considered were represented exactly by \gls{nurbs} surfaces.
For the Gaussian quadrature rule,
we relied on the well-known Gauss-Legendre one.
As the regularity of the integrand involved in the evaluation of the \gls{sir}
[\cref{eq:sir-generic-bcs:bis}]
has not been studied in depth,
we relied on a heuristic strategy consisting of selecting the number of
quadrature points in each
\( \parens{\nurbscoordu, \nurbscoordv} \)
direction of the \gls{nurbs} surface to obtain a spatial sampling rate
equivalent to the sampling rate considered for the time dimension.
By doing so,
it was observed that the resulting accuracy was not bound
by the quadrature rule in the cases studied.

Because no analytic expressions could be derived for the field signals,
we relied on a very high sampling rate
to evaluate the analytic \gls{sir},
the excitation waveform,
and their discrete convolution.
We used a reference sampling rate of \SI{20.01}{\tera\hertz}
for the \num{30}-\si{\mega\hertz} case
and of \SI{20}{\tera\hertz} for the \num{80}-\si{\mega\hertz} case
to guarantee integer downsampling factors.
We also compared the \glspl{sir} obtained using
\cref{eq:sir-from-basis-sir-conv-inverse}.

\subsubsection{Spherically focused element with a rigid baffle condition}%
\label{sec:validation:analytic:spherical}

The geometry of the spherically focused transducer element considered
is defined by an active diameter \( D = 20 \wavelength \)
and a spherical radius \( R \) defined such that
\( 2R / D = 4.8 \),
namely a similar ratio to the one studied in~\cite[Fig.~4]{Arditi_UI_1981}.
The \gls{nurbs} representation and corresponding four \emph{smooth} \gls{bezier}
patches is identical to the one described in
\cref{fig:approach:nurbs:mapping-spherical}.
Our heuristic strategy to define the number of quadrature points led to
\( \parens{59 \times 91} \) Gauss-Legendre quadrature points
per \gls{bezier} patch for the sampling rate of \SI{30}{\mega\hertz},
and \( \parens{155 \times 243} \) for the sampling rate of \SI{80}{\mega\hertz}.
The three field points (A, B, and C) at which field signals were evaluated
all lie in the same plane of revolution at a depth of
\( D / 2 = 10 \wavelength \).
The lateral coordinate of the first one (A) is \( x_{\text{A}} = 0 \).
For the second one (B),
the lateral coordinate was computed such that its projection onto
the surface lies on an edge,
resulting in \( x_{\text{B}} \approx 8.1 \wavelength \).
The lateral coordinate of the last one (C) was simply set to
\( x_{\text{C}} = 2 x_{\text{B}} \)
such that its projection onto the surface lies outside.

\begin{figure*}[htb]
    \centering
    \includegraphics[scale=0.925]{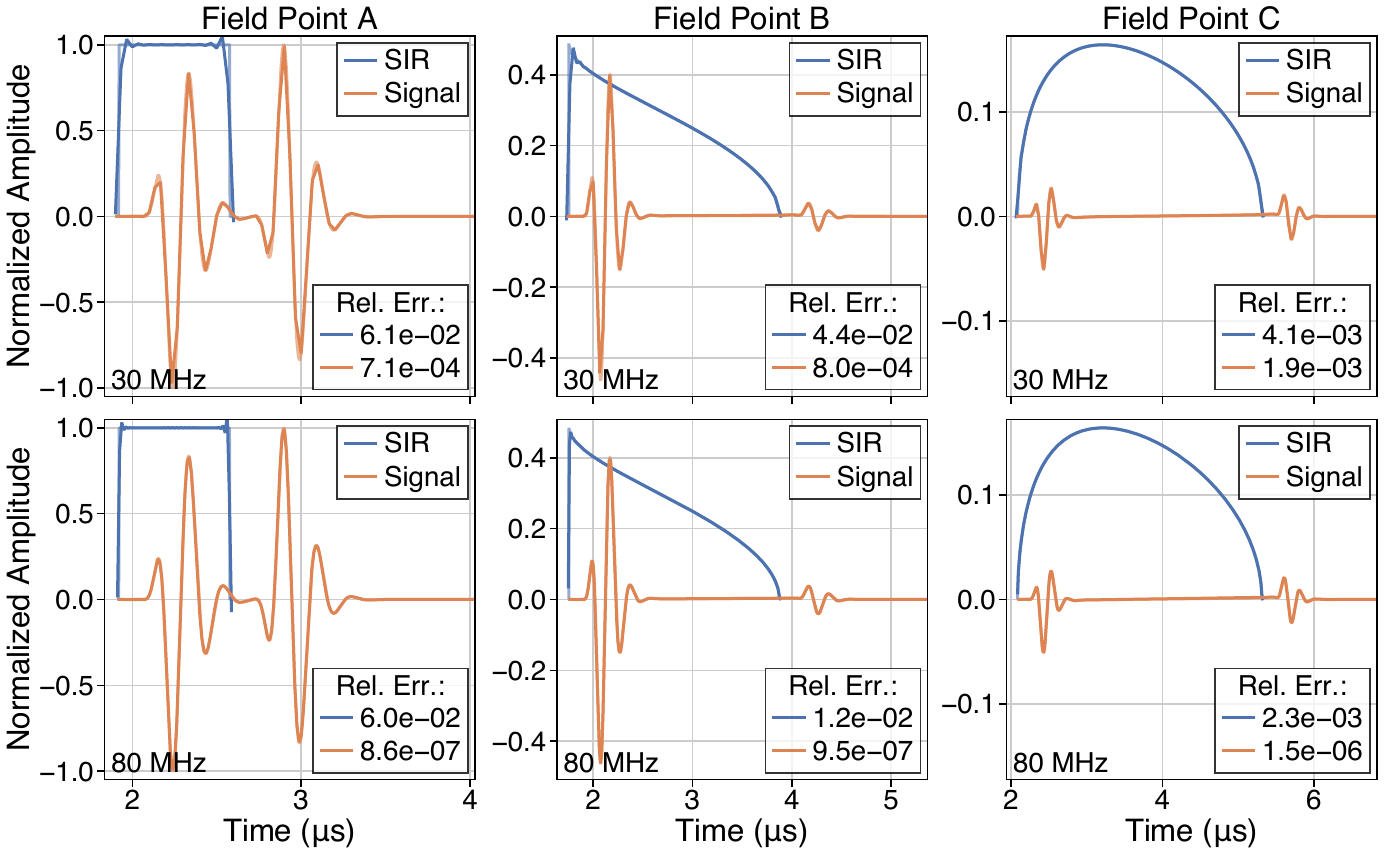}%
    \caption{%
        (Color online)
        Comparison of the \glsxtrfullpl{sir} and field signals radiated
        at different field points by
        a spherically focused transducer element with a rigid baffle condition,
        excited by a windowed-sinusoidal waveform.
        The excitation waveform is a differentiated log-normal-windowed sine wave,
        with a characteristic (center) wavelength \( \lambda \).
        The geometry of the spherical cap is defined by an active diameter
        of \( 20 \lambda \) and a spherical radius of \( 240 \lambda \).
        The three field points (A, B, C) lie in the same revolution plane
        at a depth of \( 10 \lambda \) and a lateral coordinate of
        0,
        \( 8.1 \lambda \) (\gls{ie}, projection on the edge of the surface),
        and \( 16.2 \lambda \),
        respectively.
        The proposed approach was implemented with
        a \gls{bspline} basis function of degree five,
        and was evaluated at two sampling rates of
        (first row) \SI{30}{\mega\hertz}
        and
        (second row) \SI{80}{\mega\hertz}.
        The reference \glsxtrshortpl{sir} and field signals were evaluated at a
        sampling rate of \SI{20}{\tera\hertz}.
        They are depicted with the same colors as the approximated counterparts,
        with a lower opacity.%
    }%
    \label{fig:validation:analytic:spherical}
\end{figure*}

The relative two-norm errors of the field signals at each field point
for both sampling rates and all basis functions considered are reported in
\cref{tab:validation:analytic:spherical}.
These results indicate that the basis function of the highest order
performs best (\gls{ie}, \gls{bspline}5),
with a relative error of approximately
\num{e-4} and \num{e-7}
at a sampling rate of \SI{30}{\mega\hertz} and \SI{80}{\mega\hertz},
respectively.
The field signals and \glspl{sir}
for the \gls{bspline} basis function of degree five is depicted in
\cref{fig:validation:analytic:spherical}.
Despite the very different \gls{sir} at the three field points considered,
the relative errors are similar at a given sampling rate.
One can note that the approximated \gls{sir} typically contains ripple artifacts
because the very high frequencies cannot be accounted for at such low sampling
rates.
This is especially visible for the field point A at a sampling rate of
\SI{30}{\mega\hertz}.
Yet,
the field signals do not suffer from such artifacts.
We can also observe that the approach would tend to an accurate approximation
of the \gls{sir} at much higher rates,
should such a quantity be of interest.
Field signals obtained with the other basis functions considered
can be found in
\cref{fig:app:validation:analytic:spherical:nearest,%
fig:app:validation:analytic:spherical:linear,%
fig:app:validation:analytic:spherical:keys,%
fig:app:validation:analytic:spherical:bspline3,%
fig:app:validation:analytic:spherical:omoms3,%
}.

\begin{table*}[ht]
    \caption{%
        Relative two-norm errors of field signals radiated by a
        spherically focused element with a rigid baffle condition.%
    }%
    \label{tab:validation:analytic:spherical}%
    \input{tables/simulation/tab-exp-spherical-jasa.tex}%
\end{table*}

\subsubsection{Rectangular element with a soft baffle condition}%
\label{sec:validation:analytic:rectangle-soft}

The geometry of the rectangular transducer element considered is defined by
a width of \( \wavelength \) and a height of \( 10 \wavelength \),
chosen to reflect typical width-to-height ratios of transducer elements
composing linear arrays.
The \gls{nurbs} representation is simply a bilinear surface,
resulting in a single \emph{smooth} \gls{bezier} patch
onto which the Gauss-Legendre quadrature rule was deployed.
Our heuristic strategy to define the number of quadrature points led to
\( \parens{7 \times 59} \) Gauss-Legendre quadrature points
for the sampling rate of \SI{30}{\mega\hertz},
and \( \parens{17 \times 155} \) for the sampling rate of \SI{80}{\mega\hertz}.
The three field points (A, B, and C) at which field signals were evaluated
all lie in a plane parallel to the element.
They were positioned at a depth of \( \wavelength / 2 \)
and an elevation of \( \wavelength / 2 \),
with lateral coordinates of \( 0 \), \( \wavelength / 2 \),
and \( \wavelength \), respectively.
Note that this is quite an extreme case that was selected on purpose as
this results in \glspl{sir} with very high frequencies.

\begin{figure*}[htb]
    \centering
    \includegraphics[scale=0.925]{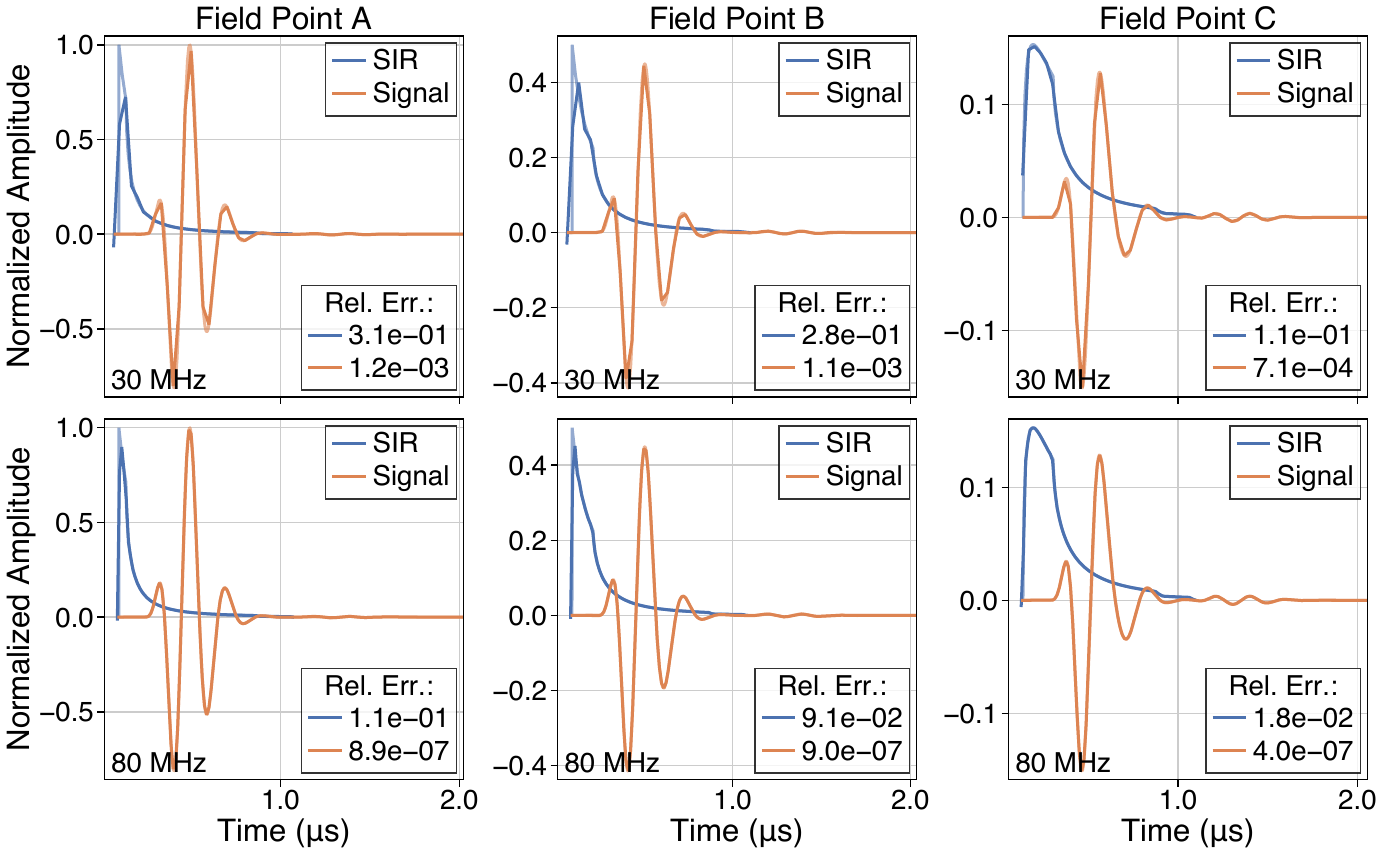}%
    \caption{%
        (Color online)
        Comparison of the \glsxtrfullpl{sir} and field signals radiated
        at different field points by
        a rectangular transducer element with a soft baffle condition,
        excited by a windowed-sinusoidal waveform.
        The excitation waveform is a differentiated log-normal-windowed sine wave,
        with a characteristic (center) wavelength \( \lambda \).
        The geometry of the rectangular element is defined by
        a width of \( \lambda \) and a height of \( 10 \lambda \).
        The three field points (A, B, C) lie in a plane parallel to the element.
        They were positioned at a depth of \( \lambda / 2 \)
        and an elevation of \( \lambda / 2 \),
        with lateral coordinates of 0, \( \lambda / 2 \), and \( \lambda \),
        respectively.
        The proposed approach was implemented with
        a \gls{bspline} basis function of degree five,
        and was evaluated at two sampling rates of
        (first row) \SI{30}{\mega\hertz}
        and
        (second row) \SI{80}{\mega\hertz}.
        The reference \glsxtrshortpl{sir} and field signals were evaluated at a
        sampling rate of \SI{20}{\tera\hertz}.
        They are depicted with the same colors as the approximated counterparts,
        with a lower opacity.%
    }%
    \label{fig:validation:analytic:rectangle-soft}
\end{figure*}

The relative two-norm errors of the field signals at each field point
for both sampling rates and all basis functions considered are reported in
\cref{tab:validation:analytic:rectangle-soft}.
As for the spherically focused cased,
the higher-order basis function performs best.
In general,
the order of relative errors are also comparable between the two element
shapes.
The field signals and \glspl{sir}
for the \gls{bspline} basis function of degree five are depicted in
\cref{fig:validation:analytic:rectangle-soft}.
Similarly to the spherically focused case,
ripples due to the high-frequency components of the \gls{sir}
can be observed in the approximations.
Field signals obtained with the other basis functions considered
can be found in
\cref{fig:app:validation:analytic:rectangle-soft:nearest,%
fig:app:validation:analytic:rectangle-soft:linear,%
fig:app:validation:analytic:rectangle-soft:keys,%
fig:app:validation:analytic:rectangle-soft:bspline3,%
fig:app:validation:analytic:rectangle-soft:omoms3,%
}.

The same experiment was also carried out for a rigid baffle condition.
All results are consistent for the two different baffle conditions.
The detailed results are reported in supplementary material
(\cref{sec:validation:analytic:rectangle-rigid}).

\begin{table*}[ht]
    \caption{%
        Relative two-norm errors of field signals radiated by a
        rectangular element with a soft baffle condition.%
    }%
    \label{tab:validation:analytic:rectangle-soft}%
    \input{tables/simulation/tab-exp-rect-soft-jasa.tex}
\end{table*}

%% file: tables/simulation/tab-exp-spherical-jasa.tex
\sisetup{
    table-format = 1.2e+1,  %
    retain-zero-exponent,
}%
\begin{tabular}{c c S S S S S S}
    \hline\hline
    {Freq.}
    & {Point}
    & {Nearest}
    & {Linear}
    & {Keys}
    & {\Gls{bspline}3}
    & {\glsxtrshort{o-moms}3}
    & {\Gls{bspline}5}
    \\
    \hline
    \multirow{3}{*}{\rotatebox[origin=c]{90}{\SI[retain-zero-exponent=false]{30}{\mega\hertz}}}
    & A
    & 2.89e-01 & 9.50e-02 & 2.14e-02 & 4.38e-03 & 1.71e-03 & 7.13e-04
    \\
    & B
    & 8.91e-02 & 9.06e-02 & 1.94e-02 & 4.64e-03 & 1.27e-03 & 7.96e-04
    \\
    & C
    & 1.41e+00 & 1.99e-01 & 1.06e-01 & 1.57e-02 & 6.96e-03 & 1.93e-03
    \\
    \hline
    \multirow{3}{*}{\rotatebox[origin=c]{90}{\SI[retain-zero-exponent=false]{80}{\mega\hertz}}}
    & A
    & 1.17e-01 & 1.40e-02 & 6.02e-04 & 6.13e-05 & 2.35e-05 & 8.62e-07
    \\
    & B
    & 1.25e-02 & 1.37e-02 & 4.51e-04 & 6.12e-05 & 2.34e-05 & 9.50e-07
    \\
    & C
    & 1.46e-01 & 1.53e-02 & 2.01e-03 & 1.03e-04 & 3.50e-05 & 1.55e-06
    \\
    \hline\hline
\end{tabular}

%% file: tables/simulation/tab-exp-rect-soft-jasa.tex
\sisetup{
    table-format = 1.2e+1,  %
    retain-zero-exponent,
}%
\begin{tabular}{c c S S S S S S}
    \hline\hline
    {Freq.}
    & {Point}
    & {Nearest}
    & {Linear}
    & {Keys}
    & {\Gls{bspline}3}
    & {\glsxtrshort{o-moms}3}
    & {\Gls{bspline}5}
    \\
    \hline
    \multirow{3}{*}{\rotatebox[origin=c]{90}{\SI[retain-zero-exponent=false]{30}{\mega\hertz}}}
    & A
    & 1.13e-01 & 1.10e-01 & 3.05e-02 & 6.18e-03 & 2.16e-03 & 1.20e-03
    \\
    & B
    & 2.66e-01 & 1.05e-01 & 2.89e-02 & 5.80e-03 & 2.04e-03 & 1.09e-03
    \\
    & C
    & 2.44e-01 & 1.08e-01 & 1.90e-02 & 5.01e-03 & 1.48e-03 & 7.12e-04
    \\
    \hline
    \multirow{3}{*}{\rotatebox[origin=c]{90}{\SI[retain-zero-exponent=false]{80}{\mega\hertz}}}
    & A
    & 1.37e-02 & 1.42e-02 & 4.13e-04 & 6.15e-05 & 2.99e-05 & 8.93e-07
    \\
    & B
    & 1.06e-02 & 1.08e-02 & 3.16e-04 & 4.31e-05 & 2.12e-05 & 9.00e-07
    \\
    & C
    & 5.02e-02 & 1.15e-02 & 4.58e-04 & 3.95e-05 & 2.07e-05 & 3.95e-07
    \\
    \hline\hline
\end{tabular}

%% file: content/discussion.tex
\section{Discussion}%
\label{sec:discussion}

The results obtained from the two experiments carried out demonstrate that
the proposed approach can attain
a high accuracy \gls{wrt} analytic reference signals.
The comparison of the relative errors obtained in all cases of the second
experiment on realistic transducer element shapes
(\cref{tab:validation:analytic:spherical,%
tab:validation:analytic:rectangle-rigid,%
tab:validation:analytic:rectangle-soft})
with those of the convergence-order study on a random stream of Dirac
delta functions
(\cref{fig:validation:interpolation:convergence-order})
shows that the latter provides an upper bound on the relative error.
As such,
the first experiment provides a robust way for selecting appropriate basis
functions to guarantee a desired accuracy,
assuming that suitable numerical quadrature is deployed.

\subsection{Benefit of high-order B-spline basis functions}%
\label{sec:discussion:high-order}

In pulse-echo \gls{us} imaging,
it is typically acceptable to consider a sampling rate that enable
preserving frequencies with a relative spectrum magnitude
of approximately \SI{-60}{\decibel}.
In the experiments carried out,
this corresponds to a sampling rate of \( \SI{\sim 35}{\mega\hertz} \)
[\cref{fig:validation:waveform:spectrum}].
One can note
(\cref{fig:validation:interpolation:convergence-order})
that using a \gls{bspline} basis function of degree five
results in a relative error of approximately \SI{-60}{\decibel},
and should therefore not induce additional errors in the simulated signals
at that sampling rate.
On the other hand,
using a nearest-neighbor basis function
would require a sampling rate of \( \SI{\sim 10}{\giga\hertz} \)
to achieve the same relative error.
This is an important observation and a major advantage of the proposed method
as large-scale simulations highly benefit from low sampling rates
(because of the many discrete convolutions).
Thus,
the use of high-order basis functions is of primary interest,
even if it requires basis functions of slightly greater supports.

\subsection{Non-uniform rational B-spline representation of surfaces}%
\label{sec:discussion:nurbs}

Our choice of representing radiating surfaces as \gls{nurbs} surfaces
was primarily motivated by the fact that they enable representing quadric
surfaces exactly,
and that most \gls{us} transducer elements are at most quadric surfaces.
Other representations could be used in the proposed approach
(\gls{eg}, analytic ones),
provided that the parametrization is \emph{smooth} such that Gaussian quadrature
can be deployed.
The fact that any \gls{nurbs} surface can be decomposed into \emph{smooth}
\gls{bezier} patches guarantees that this requirement is always met.

Also,
\gls{nurbs} representations are heavily used
in \gls{cad} software.
Hence, the design of new transducer element shapes could
be evaluated directly from these \gls{nurbs} definitions without additional
processing (such as meshing),
similarly to the concept of \gls{iga}.
The optimization of transducer element shapes could also be performed directly
from their \gls{nurbs} definitions
with direct evaluation of the field quantities of interest.

\subsection{Gaussian quadrature rules}%
\label{sec:discussion:quadrature-rules}

We only considered the well-known Gauss-Legendre quadrature rule in the
present study.
Since there is no realistic scenario in which field signals would need to
be computed onto the radiating surface
(\gls{ie}, non-singular integrand),
the Gauss-Legendre quadrature rule guarantees high-order accuracies.
Even though not reported here,
we also evaluated the performance of typical (low-order) quadrature rules
such as the midpoint, trapezoidal, and \gls{simpson} ones.
As expected,
they all performed much less efficiently than the Gauss-Legendre one.
We also conducted some preliminary evaluations using
the Gauss-Legendre-Lobatto rule.
This quadrature rule may be promising in the context of radiating surfaces
as it includes the endpoints of the integration interval
(\gls{ie}, the surface edges),
at the cost of a slightly reduced accuracy than the Gauss-Legendre rule.
However,
we could not yet conclude on which one performs generally best.
An in-depth study on the regularity of the integrand should be performed
to obtain an optimal selection of the quadrature point number.

\subsection{Comparison with other strategies}%
\label{sec:discussion:comparisons}

It is interesting to note that the proposed approach can be used to interpret
previously proposed discretized approaches,
as different surface representations,
quadrature rules,
and basis functions can be deployed.
For instance,
\citet{Piwakowski_JASA_1989} proposed
an approach to compute the \gls{sir} by discretizing the radiating surface
into many ideal points and to average
the Dirac delta functions between each time samples.
This strategy can be considered in the proposed formulation,
namely by using a midpoint quadrature rule
together with a nearest-neighbor basis function.
This would of course result in much worse accuracy than using a high-order
Gaussian quadrature rule together with a high-order basis function.
In general,
the strategy of averaging \gls{sir} values was proposed and used in many
approaches~\cite{Arditi_UI_1981,Piwakowski_JASA_1989,Dhooge_JASA_1997,%
Jensen_JCA_2001}.
Thus,
these approaches could also benefit from the formulation in basis functions
proposed here,
in particular to benefit from higher-order basis functions.

From a pure computational perspective,
simulation methods based on the \gls{sir} are particularly suited
for parallel implementations such as \gls{gpu}-based ones.
A critical point when it comes to \gls{gpu}-based implementations is
the complexity of arithmetic operations.
Analytic expressions for the \gls{sir} of specific radiating surfaces
involve complex operations such as hyperbolic and inverse hyperbolic operations,
with many cases depending on the relative positioning of field points
\gls{wrt} radiating surfaces.
Such operations would typically not be ideally suited for \gls{gpu}-based
implementations.
The proposed formulation relies on many more operations of much lower
arithmetic complexity,%
\footnote{The most complicated arithmetic operations are the square root
and the cosine in the case of soft baffle conditions.}
and is as such better suited for \gls{gpu}-based implementations.
So far,
we did not perform an exhaustive benchmark of our current implementation
against well-known software such as \gls{fieldii}.
Our initial experiments in the context of pulse-echo \gls{sa} imaging
indicated that we could reduce the computing time by approximately
two orders magnitude,
considering consumer-level \glspl{cpu} (for \gls{fieldii})
and \glspl{gpu} (for the proposed approach).
Note that the initial goal of developing this approach
was to enable us generating a sufficient amount of data
for the purpose of training \gls{cnn}-based image reconstruction methods
in the context of \gls{us} imaging~\cite{Perdios_ARXIV_2020a,Perdios_TMI_2021}.

%% file: content/conclusion.tex
\section{Conclusion}%
\label{sec:conclusion}

We proposed a spline-based \gls{sir} approach for the simulation
of field signals radiated by arbitrary shapes
embedded in both rigid and soft baffles
and excited by bandpass waveforms.
This approach consists of representing a transducer surface
as a \gls{nurbs} surface
and decomposing it into smooth \gls{bezier} patches onto which
high-order Gaussian quadrature rules can be deployed
to approximate the surface integral involved in the computation of the \gls{sir}.
Using high-order \gls{bspline} bases to express the \gls{sir},
the basis \gls{sir} amounts to a sum of shifted-and-weighted basis functions
that depend on the positions and weights of the quadrature points.
The resulting field signal is then obtained by the convolution
of the basis coefficients,
derived from the excitation waveform,
and the basis \gls{sir}.
The use of \gls{nurbs} enables accurate representations of complex surfaces,
and common transducer shapes
can be represented exactly.
High-order Gaussian quadrature rules enable using fewer quadrature points
to attain a desired accuracy than low-order ones commonly used
(\gls{eg}, midpoint or trapezoidal).
High-order \gls{bspline} basis functions enables using a simulation sampling
rate identical to the one required to represent the excitation waveform
accurately.
The numerical experiments demonstrated that the proposed approach
can attain errors as low as the sampling errors of the excitation waveform.
The extension to transducer arrays and pulse-echo settings is straightforward.
This approach is also well-suited to parallel implementations.
An initial \gls{gpu} implementation enabled us to reduce the computing time
by up to two orders of magnitude compared with the well-known \gls{fieldii}
simulator.

%% file: tail/acknowledgments.tex
\section*{Acknowledgments}%
\label{sec:acknowledgments}

This work was supported in part by the Swiss National Science Foundation
under Grant 205320\_175974.

%% file: content/sup-introduction.tex
\section{Introduction}%
\label{sec:sup:introduction}

\begin{figure*}[htb]
    \centering
    \includegraphics[scale=0.925]{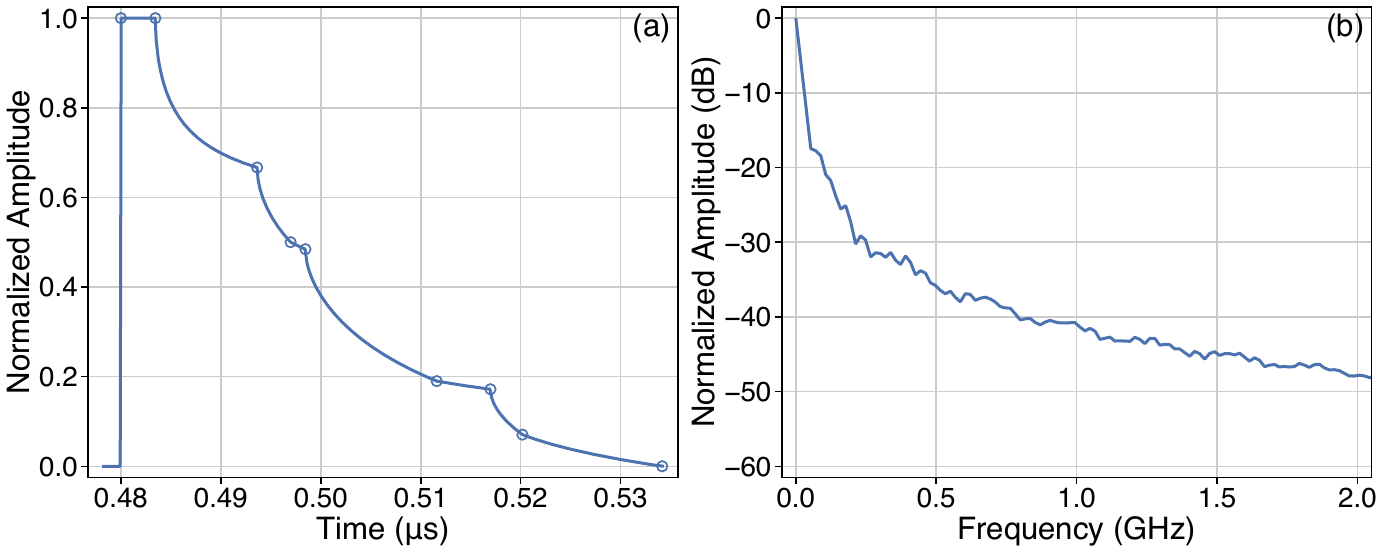}%
    \phantomsubfloatprevspace%
    \phantomsubfloat{\label{fig:introduction:sir-example:sir}}%
    \phantomsubfloat{\label{fig:introduction:sir-example:spectrum}}%
    \caption{%
        (Color online)
        Example of
        (a)
        a \glsxtrfull{sir} and
        (b)
        corresponding frequency spectrum
        for a rectangular transducer aperture.
        The rectangle has a width of \( \wavelength \) at \SI{5}{\mega\hertz}
        and a height-to-width ratio of \num{1.6}.
        The field point at which the \glsxtrshort{sir} was evaluated
        is positioned relatively to the element center at
        \( (\wavelength / \num{5}, \wavelength / \num{5}, 2.5 \wavelength) \)
        in the width, height, and depth directions.
        Characteristic abrupt slope changes are circled
        in (a).
        The first one represents the first time of arrival of a spherical shell
        centered at the position of the evaluation field point.
        The eight following ones represent time instants at which
        the spherical shell intersects the edges and vertices
        of the rectangular surface.%
    }%
    \label{fig:introduction:sir-example}
\end{figure*}

A typical example of a \gls{sir} and corresponding frequency spectrum
is shown in~\cref{fig:introduction:sir-example} for
a rectangular element designed to work at a driving frequency
of \( \SI{\sim 5}{\mega\hertz} \),
from which it is clear that a sampling rate of a few \si{\giga\hertz}
is required.

%% file: content/sup-approach.tex
\section{Proposed Approach}%
\label{sec:sup:approach}

\subsection{Non-uniform rational B-spline surface representations}%
\label{sec:sup:approach:nurbs}

\begin{figure}[htb]
    \centering
    \includegraphics[scale=0.925]{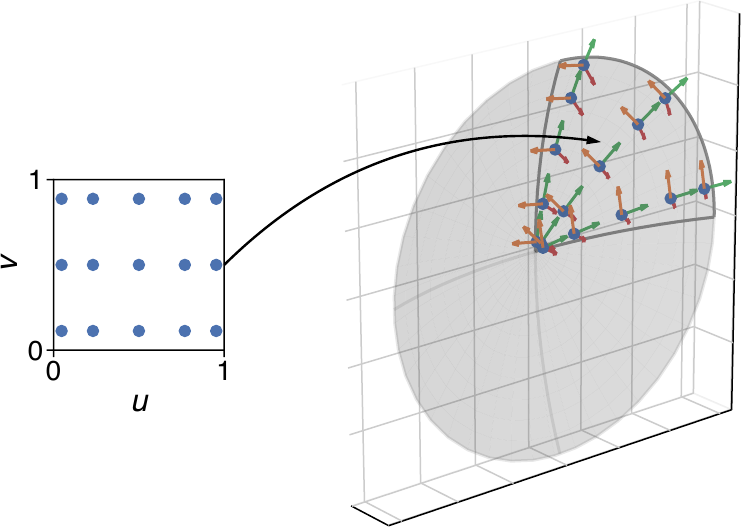}%
    \caption{%
        (Color online)
        Example of a \glsxtrfull{nurbs} surface representing a spherical cap,
        decomposed into four smooth rational biquadratic \gls{bezier} patches.
        The \glsxtrshort{nurbs} mapping from the parametric space to the physical space
        is illustrated with a parent element comprising (5\,×\,3) Gauss-Legendre
        quadrature points.
        The two (unit) tangent vectors and the (unit) normal vector
        at each quadrature point in the physical space are depicted by green,
        orange and red arrows, respectively.%
    }%
    \label{fig:approach:nurbs:mapping-spherical}
\end{figure}

As most transducer elements composing conventional \gls{us} transducers
are (at most) quadric surfaces,
the assumption that the radiating surface can be represented exactly by
a \gls{nurbs} surface is generally valid in \gls{us} imaging
since \gls{nurbs} surfaces can represent quadric surfaces
exactly~\cite[Chap.~8]{Piegl_BOOK_1997}.
For instance,
transducer elements composing linear and phased arrays are cylindrical shells
(\gls{ie}, elliptic cylinders),
those composing convex arrays are toroidal shells
(\gls{ie}, hyperbolic paraboloids),
and those composing \ndim{2} matrix arrays are simple rectangular planes.
Spherically focused transducer elements are spherical caps.
An example of the latter is shown
in \cref{fig:approach:nurbs:mapping-spherical},
for which the \gls{nurbs} surface definition can be obtained by revolving
a circular arc (\gls{ie}, represented by \gls{nurbs} curve).
The \gls{nurbs} surface can then be decomposed into four rational
biquadratic \gls{bezier} smooth patches.
Each of these smooth patches is defined by \cref{eq:nurbs-surface-mapping}
onto the characteristic unit square space
\( \bracks{\nurbsintervala, \nurbsintervalb}^{2} \).
As they are smooth surfaces,
Gaussian quadrature rules can be deployed onto the \ndim{2} parametric space
(also referred to as parent element) by means of the tensor product.
Resulting quadrature coordinates can then be mapped onto the radiating surface
exactly using \cref{eq:nurbs-surface-mapping}.
This process is also illustrated in
\cref{fig:approach:nurbs:mapping-spherical},
in which the Gauss-Legendre quadrature rule was considered to obtain
the quadrature coordinates of a \( (5 \times 3) \) grid.
Note that Gaussian quadrature rules are typically defined on a characteristic
\( \bracks{-1, 1} \) interval~\cite[Sec.~25.4]{AbramowitzAndStegun_BOOK_1964},
but can be mapped to arbitrary intervals.

\begin{figure}[htb]
    \centering
    \includegraphics[scale=0.925]{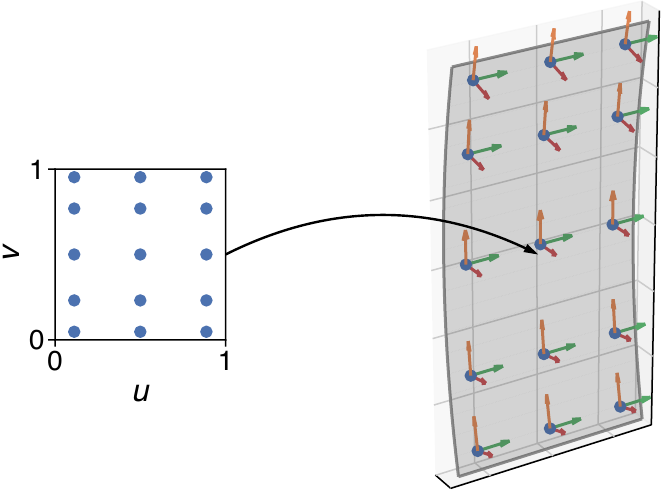}%
    \caption{%
        (Color online)
        Example of a \glsxtrfull{nurbs} surface representing a cylindrical shell
        that is also a smooth rational \gls{bezier} of degree (1, 2),
        that is,
        linear in \( u \) and quadratic in \( v \)
        directions.
        The \glsxtrshort{nurbs} mapping from the parametric space to the physical space
        is illustrated with a parent element comprising (3×5) Gauss-Legendre
        quadrature points.
        The two (unit) tangent vectors and the (unit) normal vector
        at each quadrature point in the physical space are depicted by green,
        orange and red arrows, respectively.%
    }%
    \label{fig:approach:nurbs:mapping-cylindrical}
\end{figure}

Thanks to their small (or nonexistent) curvatures,%
all other aforementioned transducer element shapes can be represented by
\gls{nurbs} surfaces that are themselves single rational \gls{bezier} patches
already,
and thus do not require any further decomposition.
An example of such as case for the cylindrical shell is shown in
\cref{fig:approach:nurbs:mapping-cylindrical}.

\subsection{Complete process and implementation details}%
\label{sec:sup:approach:implementation}

\begin{figure*}[htb]
    \centering
    \includegraphics[scale=0.925]{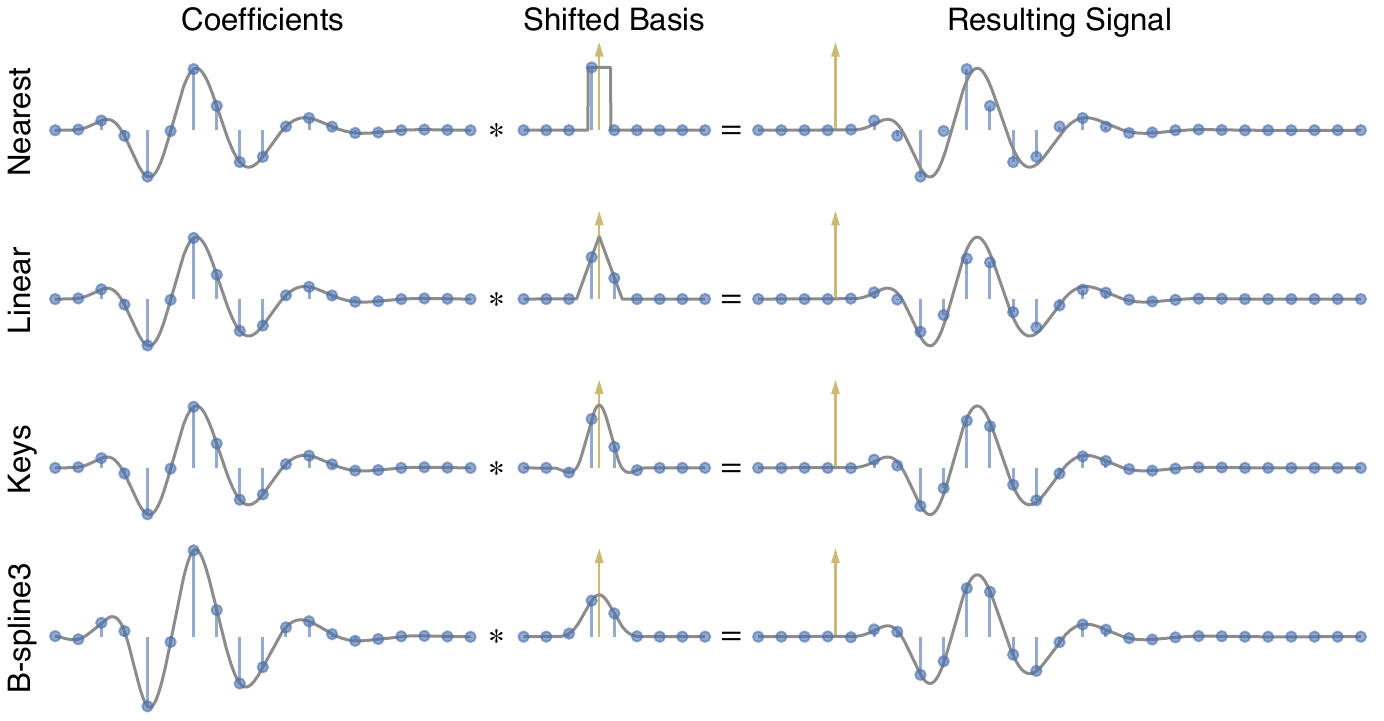}%
    \caption{%
        (Color online)
        Illustration of the generalized interpolation strategy deployed for
        approximating the signal radiated by a transducer element.
        This strategy involves the convolution of coefficients,
        evaluated from the excitation waveform,
        with shifted basis functions.
        Three interpolating basis functions are depicted,
        namely nearest-neighbor,
        linear,
        and
        (quadratic) Keys.
        One non-interpolating basis function is also shown,
        namely cubic \gls{bspline} (B-spline3).%
    }%
    \label{fig:approach:interpolation-strategy}
\end{figure*}

An illustration of the convolution involved in the evaluation of
\cref{eq:sir-generic-bandpass-signal-quadrature-discrete-geninterp-fast}
for each quadrature point is depicted in
\cref{fig:approach:interpolation-strategy}
for four different basis functions,
namely nearest-neighbor,
linear,
quadratic \gls{keys},
and cubic \gls{bspline}.
As the first three basis functions are interpolating,
their corresponding (pre-filtered) coefficients are equal to the excitation
waveform samples.
The coefficients of the cubic \gls{bspline} have larger amplitudes than
the excitation waveform to compensate for its non-interpolating property.
One can already note the major issue associated with near-neighbor and linear
interpolation,
namely a shift error and an underestimated response,
respectively.
The advantage of using a cubic \gls{bspline} over a quadratic Keys,
both having a support of four samples,
can already be noticed.
This will be confirmed by the numerical validation of the expected
convergence orders (\cref{sec:validation:interpolation}).

\clearpage
\input{content/approach-array.tex}

%% file: content/approach-array.tex
\subsection{Extension to arrays}%
\label{sec:approach:arrays}

The extension to \gls{us} transducers composed of multiple transducer elements,
typically arranged as arrays,
is straightforward.
Let us consider a generic definition of a transducer array composed of
a set of \( \sirarrayelemnb \) transducer elements.
Each element can be represented by a \gls{nurbs} surface
and decomposed into a set of smooth \gls{bezier} patches
onto which a Gaussian quadrature rule can be defined.
These elements are generally of the same shape,
which would result in similar \gls{nurbs} representations,
although this is not required.
Each \( \sirarrayelemindex \)-th element can be driven by a different
excitation waveform
\(
    \setroster{
        \sirarrayelemwaveform\parens{\sirexcitationtimeindex\sirdeltat}
    }
\),
from which the corresponding basis coefficients
\(
    \setroster{
        \sirarrayelemcoeffs\parens{\sirexcitationtimeindex\sirdeltat}
    }
\)
needs only to be evaluated once [\cref{eq:interp-pre-filtering}].
Using the principle of superposition (linear acoustics)
and \cref{eq:sir-generic-bandpass-signal-quadrature-discrete-geninterp-fast-compact},
the field signal at any field point \( \fieldpoint \) radiated by such an array
can be expressed as
\begin{align}
    \sirgenericsignal\parens{\fieldpoint, \timevar}
    &\approx
    \sum\limits_{\sirarrayelemindex=1}^{\sirarrayelemnb}
    \sum\limits_{\sirexcitationtimeindex \in \integernumbers}
    \sirarrayelemcoeffs\parens{\sirexcitationtimeindex\sirdeltat}
    \sirarrayelembsir\parens*{
        \fieldpoint,
        \timevar - \sirexcitationtimeindex\sirdeltat
    }
    \label{eq:sir-array-generic}
    ,
\end{align}
where each element basis \gls{sir} \( \sirarrayelembsir \)
can be computed using
\cref{eq:basis-sir-general}.
It is common to define explicitly the delays
\( \setroster{\sirarrayelemdelay} \)
and apodization weights
\( \setroster{\sirarrayelemapod} \)
applied across the array elements,
rather than implicitly including them in the different excitation waveforms.
Such delays and apodization weights are typically used to shape the transmit
beam (beamforming),
for instance to focus at a desired position or to steer an unfocused wavefront.
Thus,
and without loss of generality,
\cref{eq:sir-array-generic} can be rewritten as
\begin{align}
    \sirgenericsignal\parens{\fieldpoint, \timevar}
    &\approx
    \sum\limits_{\sirarrayelemindex=1}^{\sirarrayelemnb}
    \sirarrayelemapod
    \sum\limits_{\sirexcitationtimeindex \in \integernumbers}
    \sirarrayelemcoeffs\parens{\sirexcitationtimeindex\sirdeltat}
    \sirarrayelembsir\parens*{
        \fieldpoint,
        \timevar - \sirarrayelemdelay - \sirexcitationtimeindex\sirdeltat
    }
    \label{eq:sir-array-generic-delay-apod}
    .
\end{align}

In pulse-echo imaging,
it is conventionally assumed that all (identical) transducer elements composing
the array have the same electromechanical impulse response.
An identical electric excitation is also typically used,
such that the characteristic excitation waveform is identical.
Only the delays and apodization weights applied to the elements
may therefore differ such that \cref{eq:sir-array-generic-delay-apod}
can be simplified as
\begin{align}
    \sirgenericsignal\parens{\fieldpoint, \timevar}
    &\approx
    \sum\limits_{\sirexcitationtimeindex \in \integernumbers}
    \sirgenericexcitationcoeffs\parens{\sirexcitationtimeindex\sirdeltat}
    \sum\limits_{\sirarrayelemindex=1}^{\sirarrayelemnb}
    \sirarrayelemapod
    \sirarrayelembsir\parens*{
        \fieldpoint,
        \timevar - \sirarrayelemdelay - \sirexcitationtimeindex\sirdeltat
    }
    \label{eq:sir-array-fixed-waveform-delay-apod}
    .
\end{align}
Finally,
the echo signal scattered by an ideal reflector at \( \fieldpoint \)
can be easily obtained by the convolution of the field signal
of the array [\cref{eq:sir-array-fixed-waveform-delay-apod}]
with the field signal of the receive element
[\cref{eq:sir-generic-bandpass-signal-quadrature-discrete-geninterp-fast-compact}],
multiplied by the scattering amplitude~\cite{Jensen_JASA_1991}.
Note that the excitation waveform on transmit is generally different from
the one on receive,
because the first contains both the electric excitation and the
electromechanical impulse response,
whereas the second contains the electromechanical impulse response only.

%% file: content/sup-validation.tex
\section{Experiments and Results}%
\label{sec:sup:validation}

\begin{figure*}[htb]
    \centering
    \includegraphics[scale=0.925]{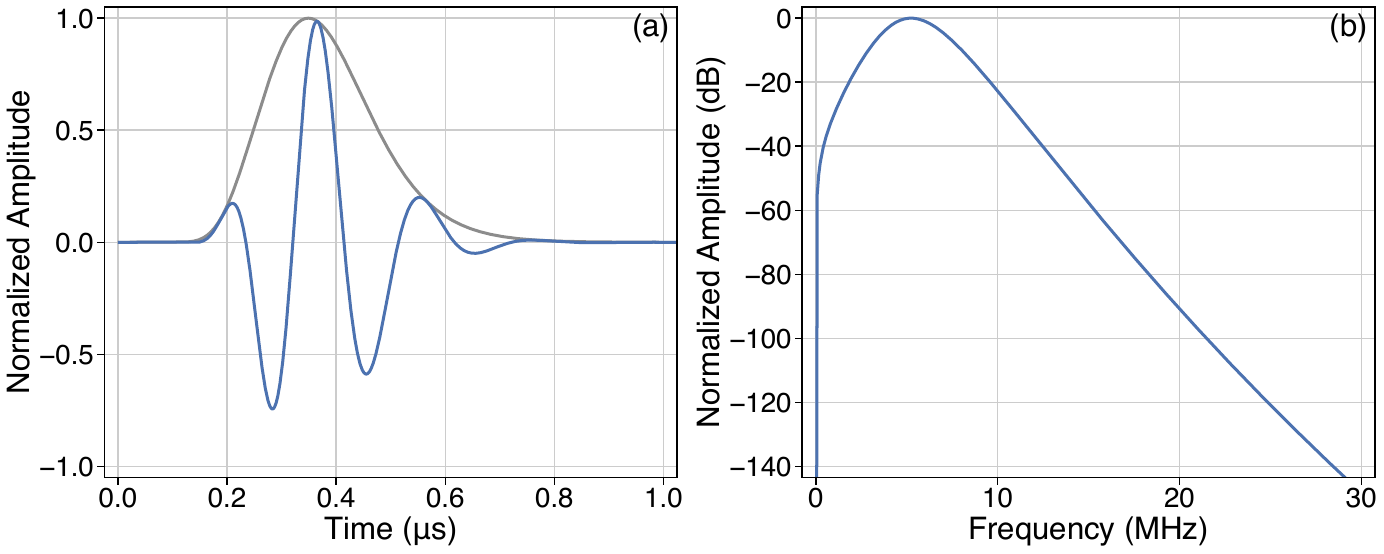}%
    \phantomsubfloatprevspace%
    \phantomsubfloat{\label{fig:validation:waveform:shape}}%
    \phantomsubfloat{\label{fig:validation:waveform:spectrum}}%
    \caption{%
        (Color online)
        Normalized amplitude of
        (a)
        the excitation waveform
        and
        (b)
        corresponding frequency spectrum
        of the pulse model considered for the numerical experiments.%
    }%
    \label{fig:validation:waveform}
\end{figure*}

The analytic expression of the time derivative
of a log-normal-modulated sinusoidal \gls{rf} pulse considered
for the experiments
(\cref{sec:validation:interpolation,sec:validation:analytic})
is shown in \cref{fig:validation:waveform}.
It is a fairly good model for the electromechanical impulse response
of transducer elements.
Note that the time derivative guarantees a zero \gls{dc} component
[\cref{fig:validation:waveform:spectrum}],
a physical property of such electromechanical impulse response.

\subsection{Validation against analytic solutions}%
\label{sec:sup:validation:analytic}

\subsubsection{Rectangular element with a rigid baffle condition}%
\label{sec:validation:analytic:rectangle-rigid}

\begin{figure*}[htb]
    \centering
    \includegraphics[scale=0.925]{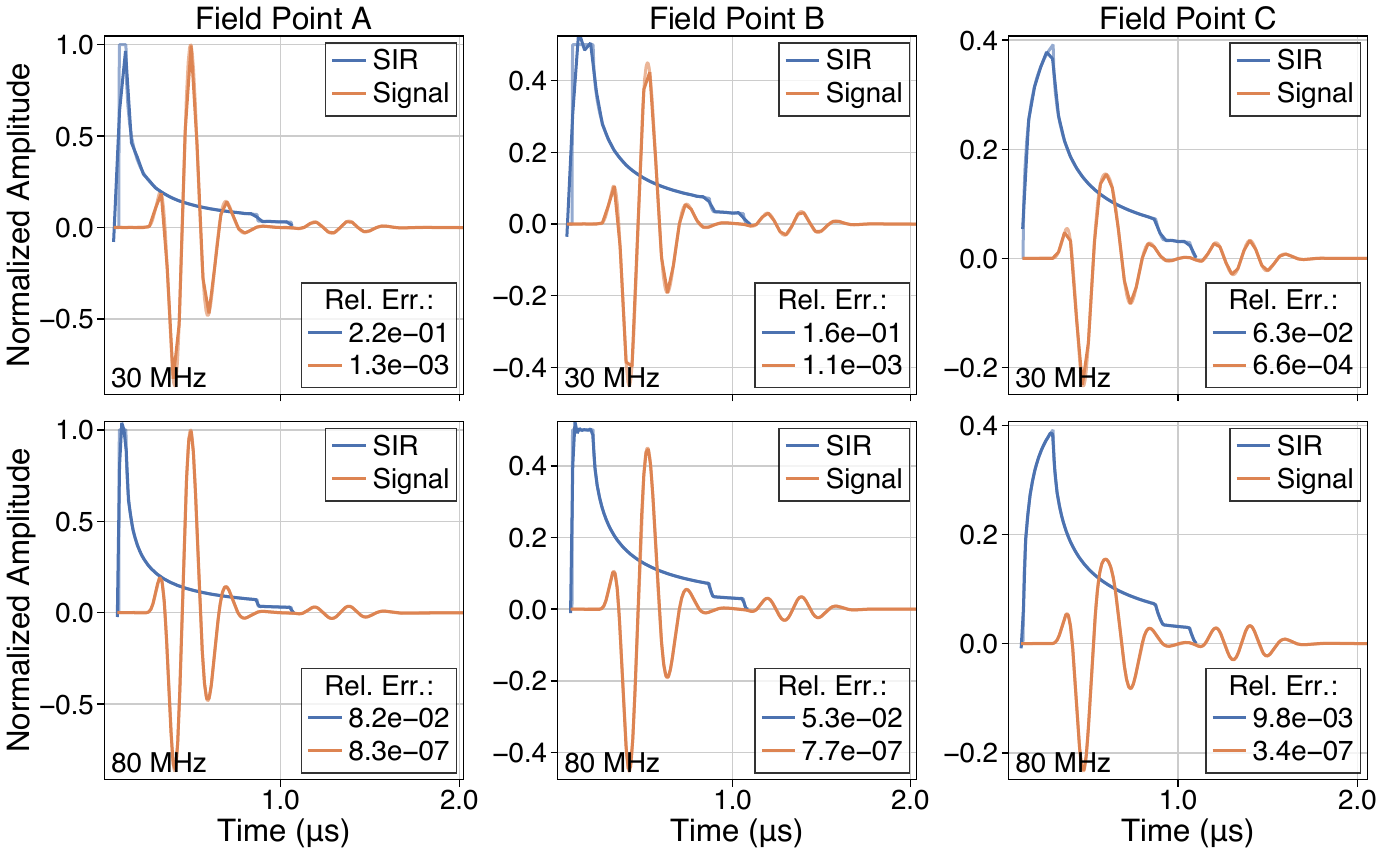}%
    \caption{%
        (Color online)
        Comparison of the \glsxtrfullpl{sir} and field signals radiated
        at different field points by
        a rectangular transducer element with a rigid baffle condition,
        excited by a windowed-sinusoidal waveform.
        The excitation waveform is a differentiated log-normal-windowed sine wave,
        with a characteristic (center) wavelength \( \lambda \).
        The geometry of the rectangular element is defined by
        a width of \( \lambda \) and a height of \( 10 \lambda \).
        The three field points (A, B, C) lie in a plane parallel to the element.
        They were positioned at a depth of \( \lambda / 2 \)
        and an elevation of \( \lambda / 2 \),
        with lateral coordinates of 0, \( \lambda / 2 \), and \( \lambda \),
        respectively.
        The proposed approach was implemented with
        a \gls{bspline} basis function of degree five,
        and was evaluated at two sampling rates of
        (first row) \SI{30}{\mega\hertz}
        and
        (second row) \SI{80}{\mega\hertz}.
        The reference \glsxtrshortpl{sir} and field signals were evaluated at a
        sampling rate of \SI{20}{\tera\hertz}.
        They are depicted with the same colors as the approximated counterparts,
        with a lower opacity.%
    }%
    \label{fig:validation:analytic:rectangle-rigid}
\end{figure*}

The same experiment as described in
\cref{sec:validation:analytic:rectangle-soft}
was carried out for the rectangular transducer element
with a rigid baffle condition.
The relative two-norm errors of the field signals at each field point
for both sampling rates and all basis functions considered are reported in
\cref{tab:validation:analytic:rectangle-rigid}.
All results are consistent for the two different baffle conditions,
even at field point C where the \glspl{sir} of the rigid and soft baffle
conditions differ greatly.
The field signals and \glspl{sir}
for the \gls{bspline} basis function of degree five are depicted in
\cref{fig:validation:analytic:rectangle-rigid}.
By comparing
\cref{fig:validation:analytic:rectangle-soft}
and
\cref{fig:validation:analytic:rectangle-rigid},
one can see that
the effect of the soft baffle condition is visible for all field points,
and in particular at field point C.
Field signals obtained with the other basis functions considered
can be found in
\cref{fig:app:validation:analytic:rectangle-rigid:nearest,%
fig:app:validation:analytic:rectangle-rigid:linear,%
fig:app:validation:analytic:rectangle-rigid:keys,%
fig:app:validation:analytic:rectangle-rigid:bspline3,%
fig:app:validation:analytic:rectangle-rigid:omoms3,%
}.

\begin{table*}[ht]
    \caption{%
        Relative two-norm errors of field signals radiated by a
        rectangular element with a rigid baffle condition.%
    }%
    \label{tab:validation:analytic:rectangle-rigid}%
    \input{tables/simulation/tab-exp-rect-rigid-jasa.tex}
\end{table*}

\clearpage
\subsubsection{Spherically focused element with a rigid baffle condition: Additional results}%
\label{sec:sup:validation:analytic:spherical}

\begin{figure*}[htb]
    \centering
    \includegraphics[scale=0.925]{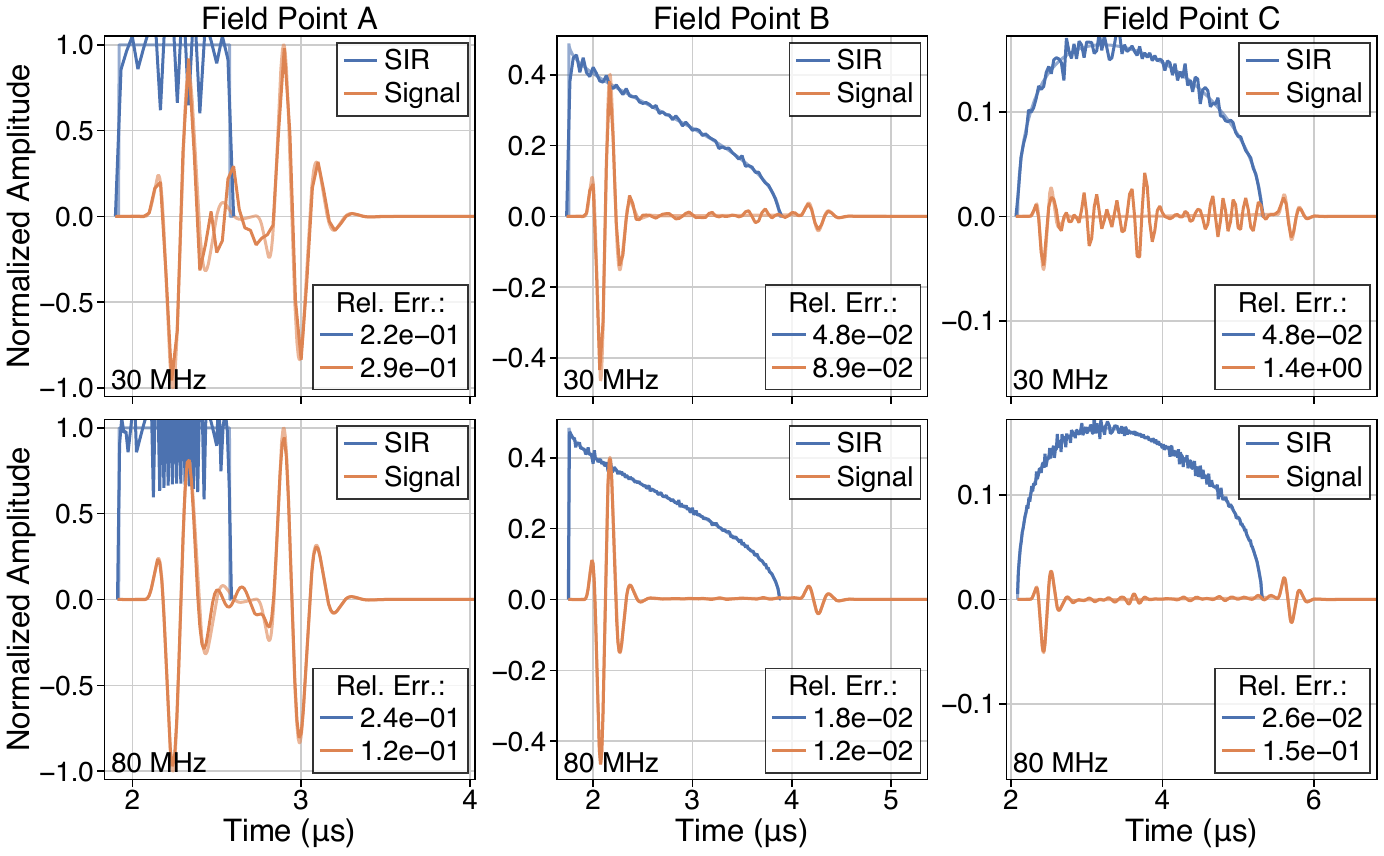}%
    \caption{%
        (Color online)
        Comparison of the \glsxtrfullpl{sir} and field signals radiated
        at different field points by
        a spherically focused transducer element with a rigid baffle condition,
        excited by a windowed-sinusoidal waveform.
        The excitation waveform is a differentiated log-normal-windowed sine wave,
        with a characteristic (center) wavelength \( \lambda \).
        The geometry of the spherical cap is defined by an active diameter
        of \( 20 \lambda \) and a spherical radius of \( 240 \lambda \).
        The three field points (A, B, C) lie in the same revolution plane
        at a depth of \( 10 \lambda \) and a lateral coordinate of
        0,
        \( 8.1 \lambda \) (\gls{ie}, projection on the edge of the surface),
        and \( 16.2 \lambda \),
        respectively.
        The proposed approach was implemented with
        a nearest-neighbor basis function (degree zero),
        and was evaluated at two sampling rates of
        (first row) \SI{30}{\mega\hertz}
        and
        (second row) \SI{80}{\mega\hertz}.
        The reference \glsxtrshortpl{sir} and field signals were evaluated at a
        sampling rate of \SI{20}{\tera\hertz}.
        They are depicted with the same colors as the approximated counterparts,
        with a lower opacity.%
    }%
    \label{fig:app:validation:analytic:spherical:nearest}
\end{figure*}

\begin{figure*}[htb]
    \centering
    \includegraphics[scale=0.925]{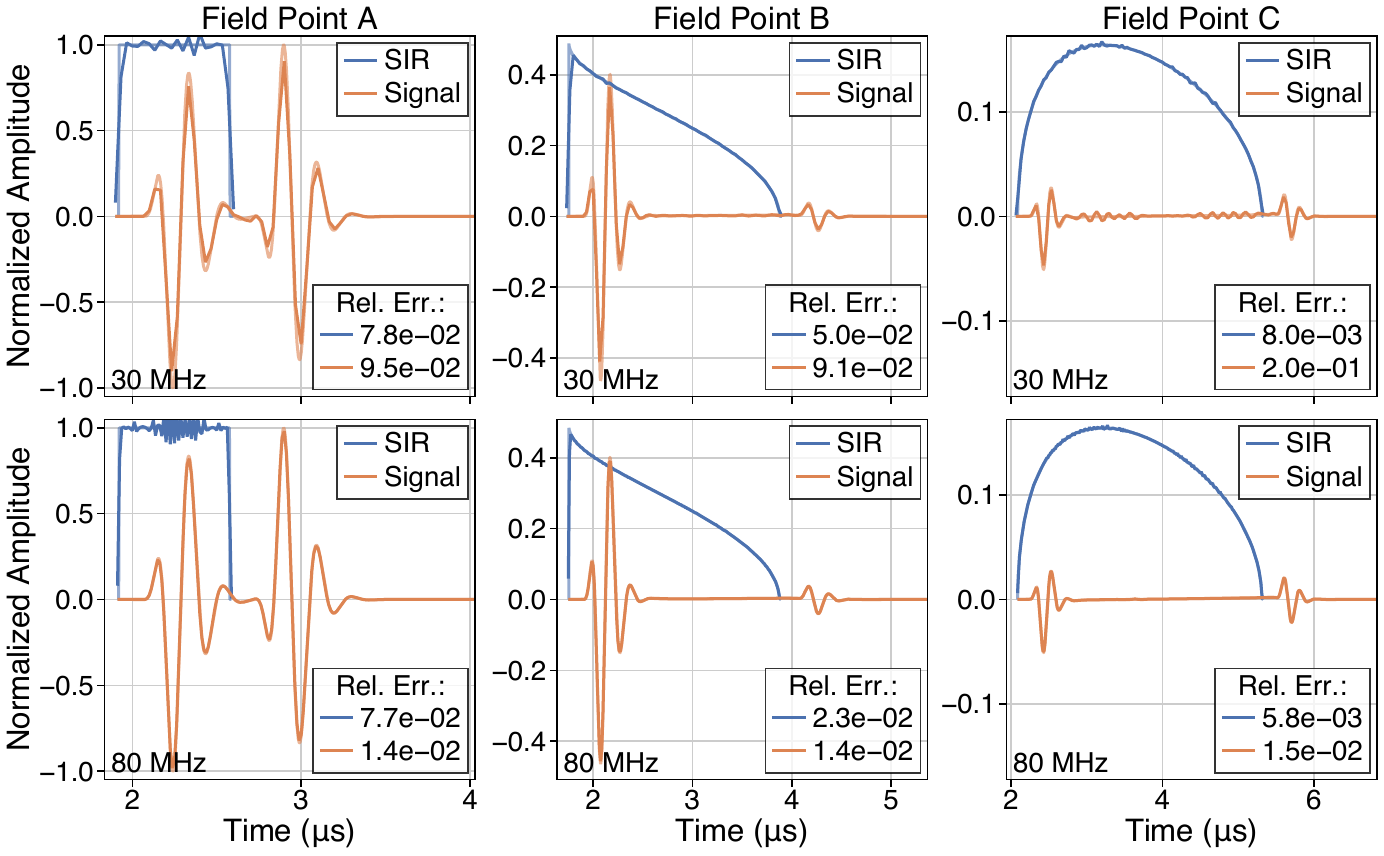}%
    \caption{(Color online) \Gls{ibid} for a linear basis function (degree one).}%
    \label{fig:app:validation:analytic:spherical:linear}
\end{figure*}

\begin{figure*}[htb]
    \centering
    \includegraphics[scale=0.925]{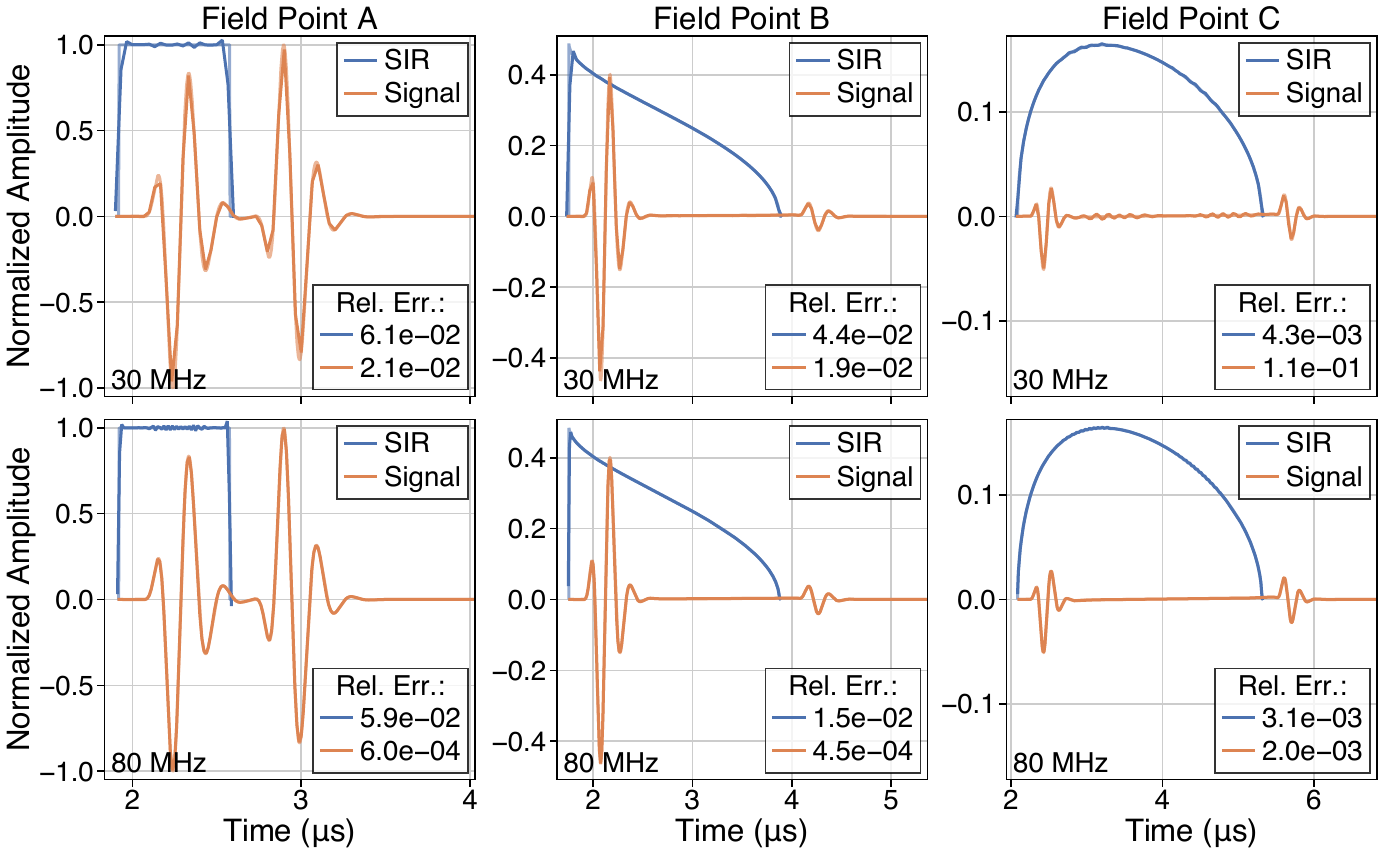}%
    \caption{(Color online) \Gls{ibid} for a quadratic \gls{keys} basis function.}%
    \label{fig:app:validation:analytic:spherical:keys}
\end{figure*}

\begin{figure*}[htb]
    \centering
    \includegraphics[scale=0.925]{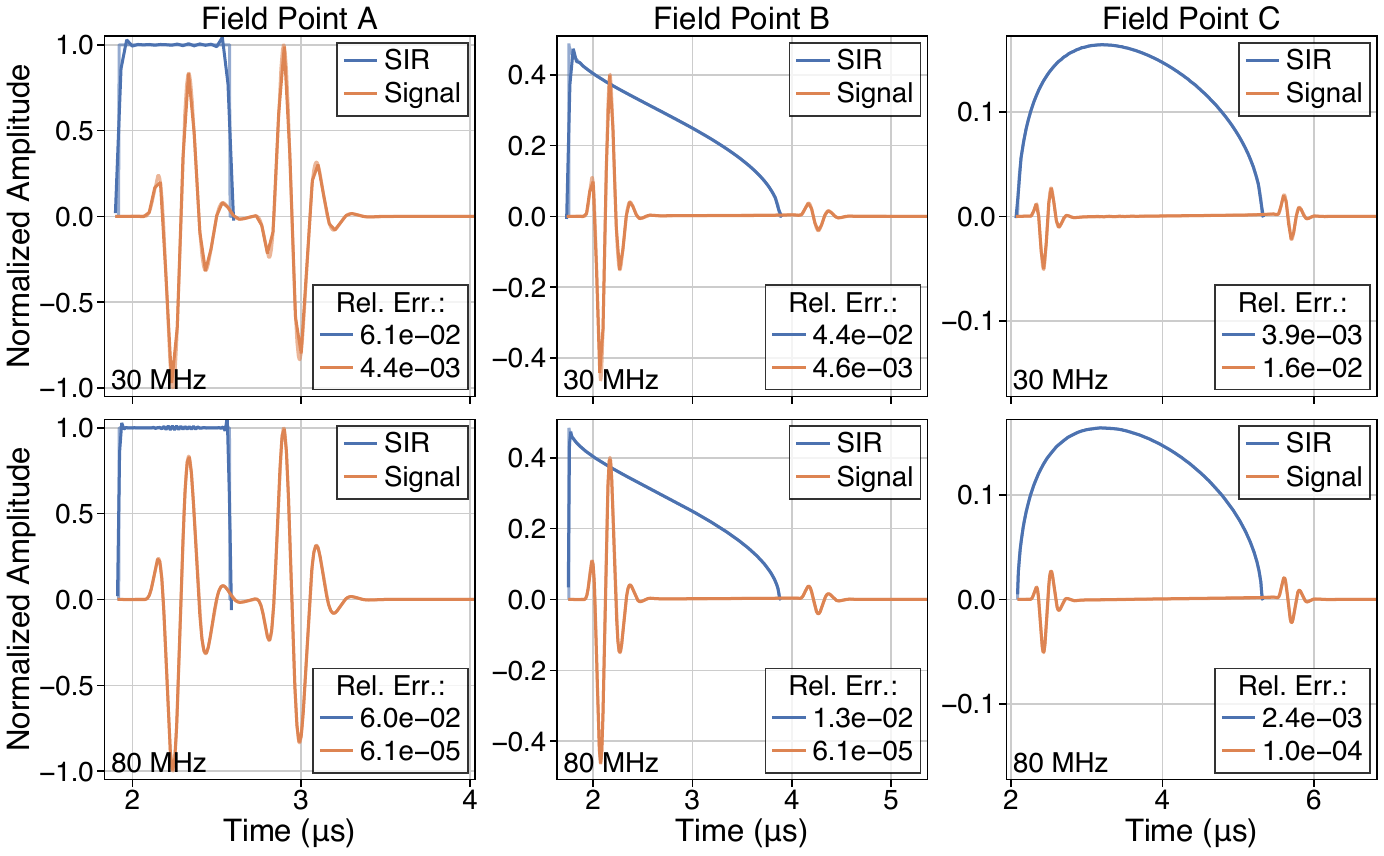}%
    \caption{(Color online) \Gls{ibid} for a cubic \gls{bspline} basis function.}%
    \label{fig:app:validation:analytic:spherical:bspline3}
\end{figure*}

\begin{figure*}[htb]
    \centering
    \includegraphics[scale=0.925]{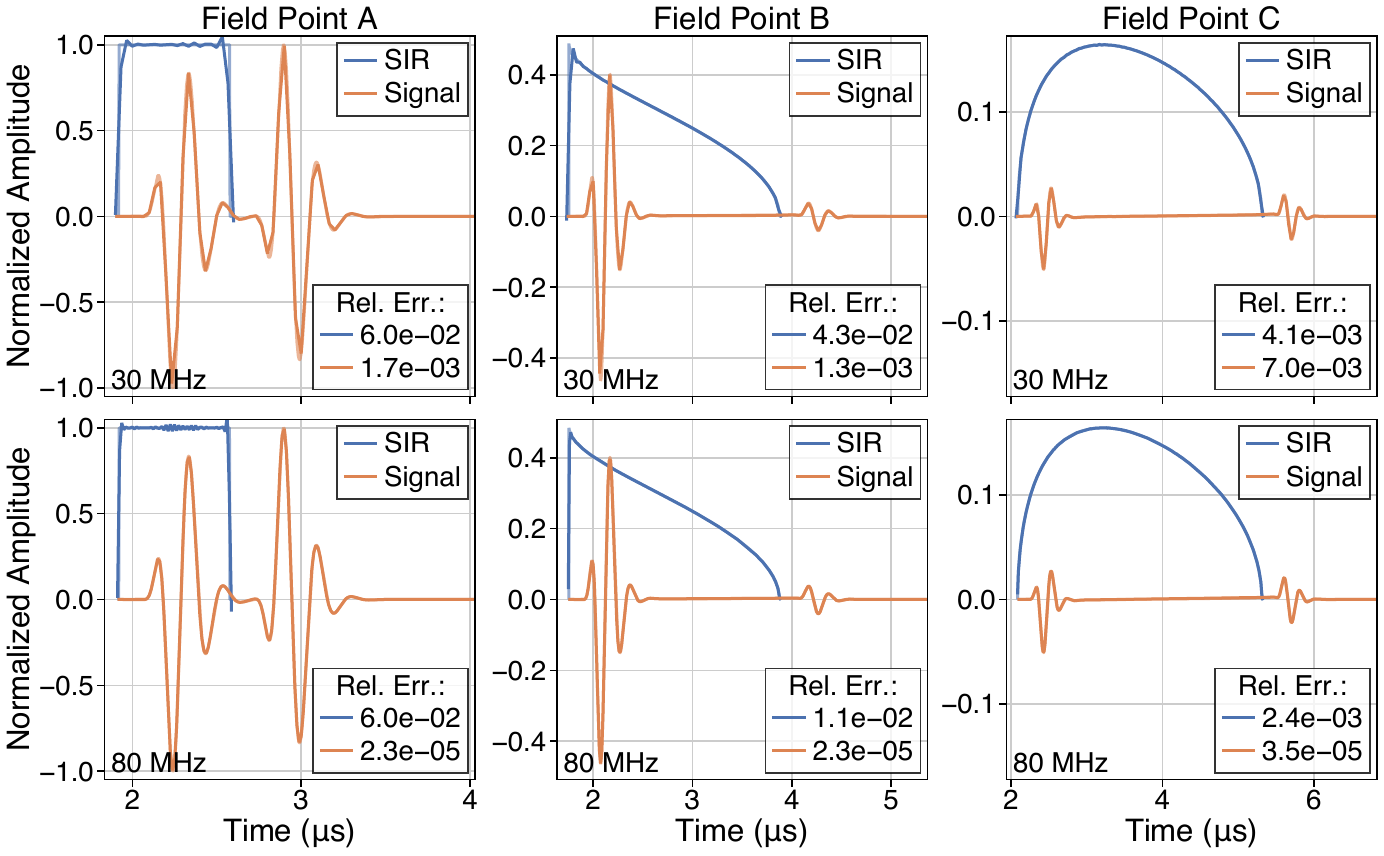}%
    \caption{(Color online) \Gls{ibid} for a cubic \glsxtrshort{o-moms} basis function.}%
    \label{fig:app:validation:analytic:spherical:omoms3}
\end{figure*}

\clearpage
\subsubsection{Rectangular element with a soft baffle condition: Additional results}%
\label{sec:sup:validation:analytic:rectangle-soft}

\begin{figure*}[htb]
    \centering
    \includegraphics[scale=0.925]{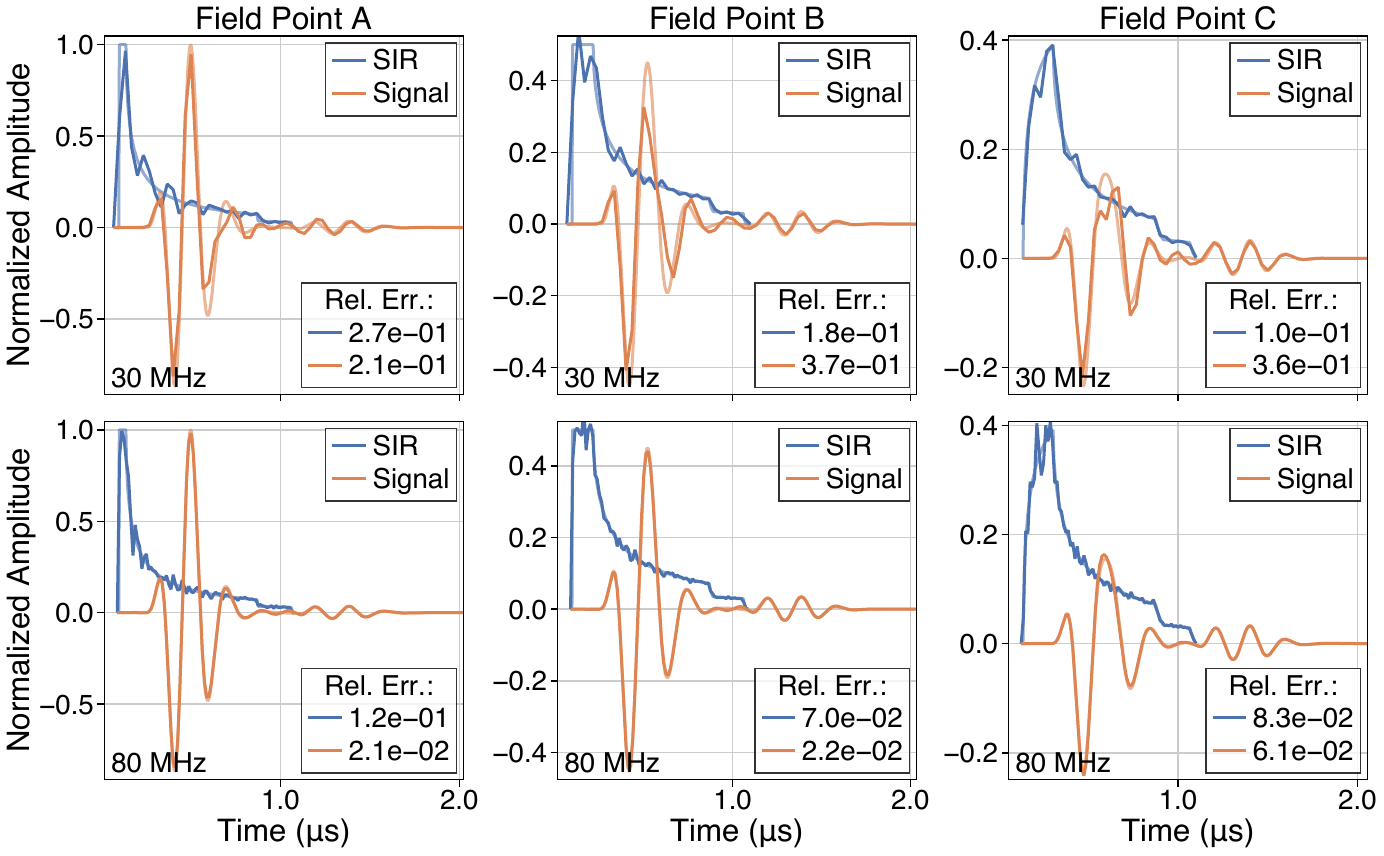}%
    \caption{%
        (Color online)
        Comparison of the \glsxtrfullpl{sir} and field signals radiated
        at different field points by
        a rectangular transducer element with a soft baffle condition,
        excited by a windowed-sinusoidal waveform.
        The excitation waveform is a differentiated log-normal-windowed sine wave,
        with a characteristic (center) wavelength \( \lambda \).
        The geometry of the rectangular element is defined by
        a width of \( \lambda \) and a height of \( 10 \lambda \).
        The three field points (A, B, C) lie in a plane parallel to the element.
        They were positioned at a depth of \( \lambda / 2 \)
        and an elevation of \( \lambda / 2 \),
        with lateral coordinates of 0, \( \lambda / 2 \), and \( \lambda \),
        respectively.
        The proposed approach was implemented with
        a nearest-neighbor basis function (degree zero),
        and was evaluated at two sampling rates of
        (first row) \SI{30}{\mega\hertz}
        and
        (second row) \SI{80}{\mega\hertz}.
        The reference \glsxtrshortpl{sir} and field signals were evaluated at a
        sampling rate of \SI{20}{\tera\hertz}.
        They are depicted with the same colors as the approximated counterparts,
        with a lower opacity.%
    }%
    \label{fig:app:validation:analytic:rectangle-soft:nearest}
\end{figure*}

\begin{figure*}[htb]
    \centering
    \includegraphics[scale=0.925]{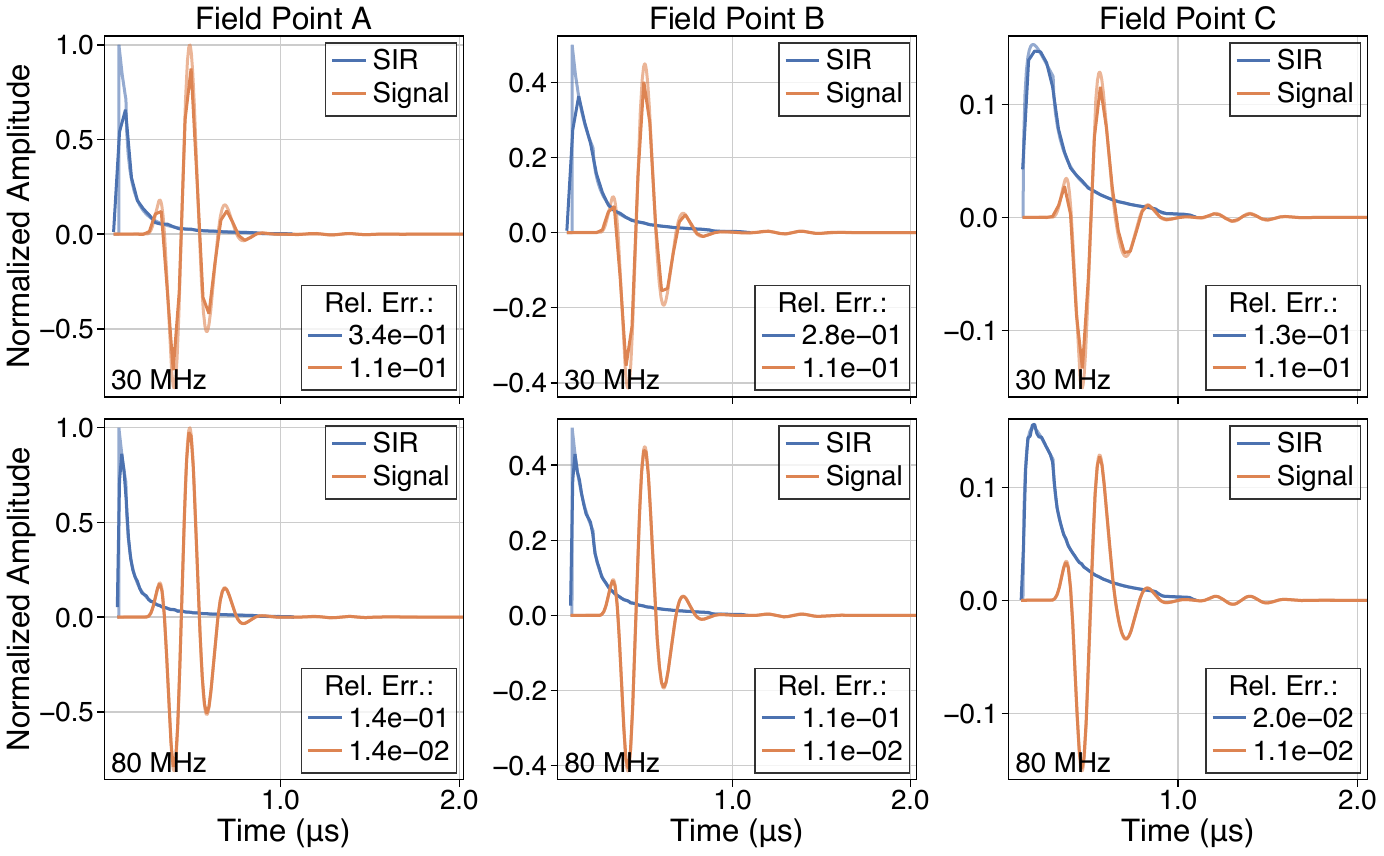}%
    \caption{(Color online) \Gls{ibid} for a linear basis function (degree one).}%
    \label{fig:app:validation:analytic:rectangle-soft:linear}
\end{figure*}

\begin{figure*}[htb]
    \centering
    \includegraphics[scale=0.925]{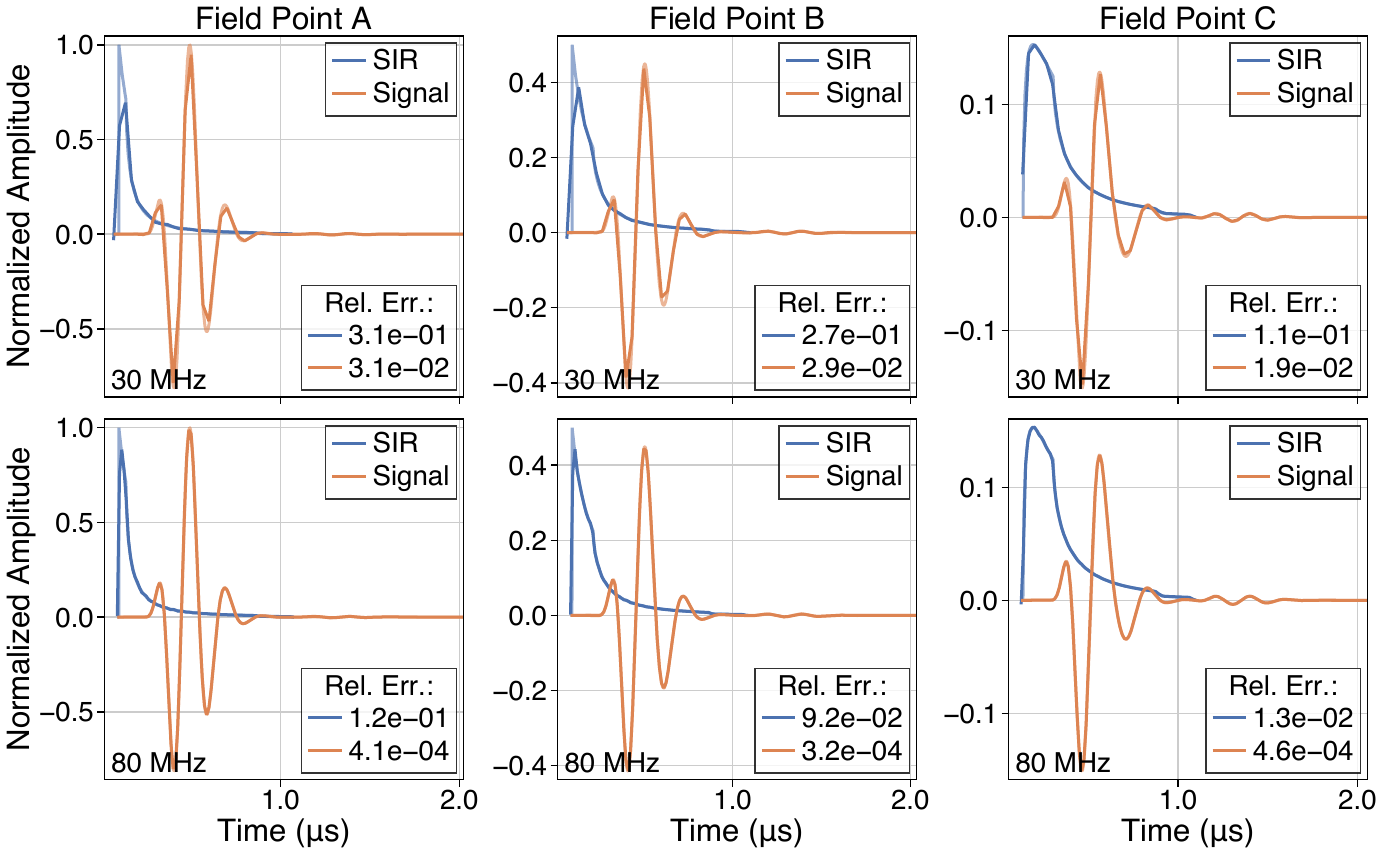}%
    \caption{(Color online) \Gls{ibid} for a quadratic \gls{keys} basis function.}%
    \label{fig:app:validation:analytic:rectangle-soft:keys}
\end{figure*}

\begin{figure*}[htb]
    \centering
    \includegraphics[scale=0.925]{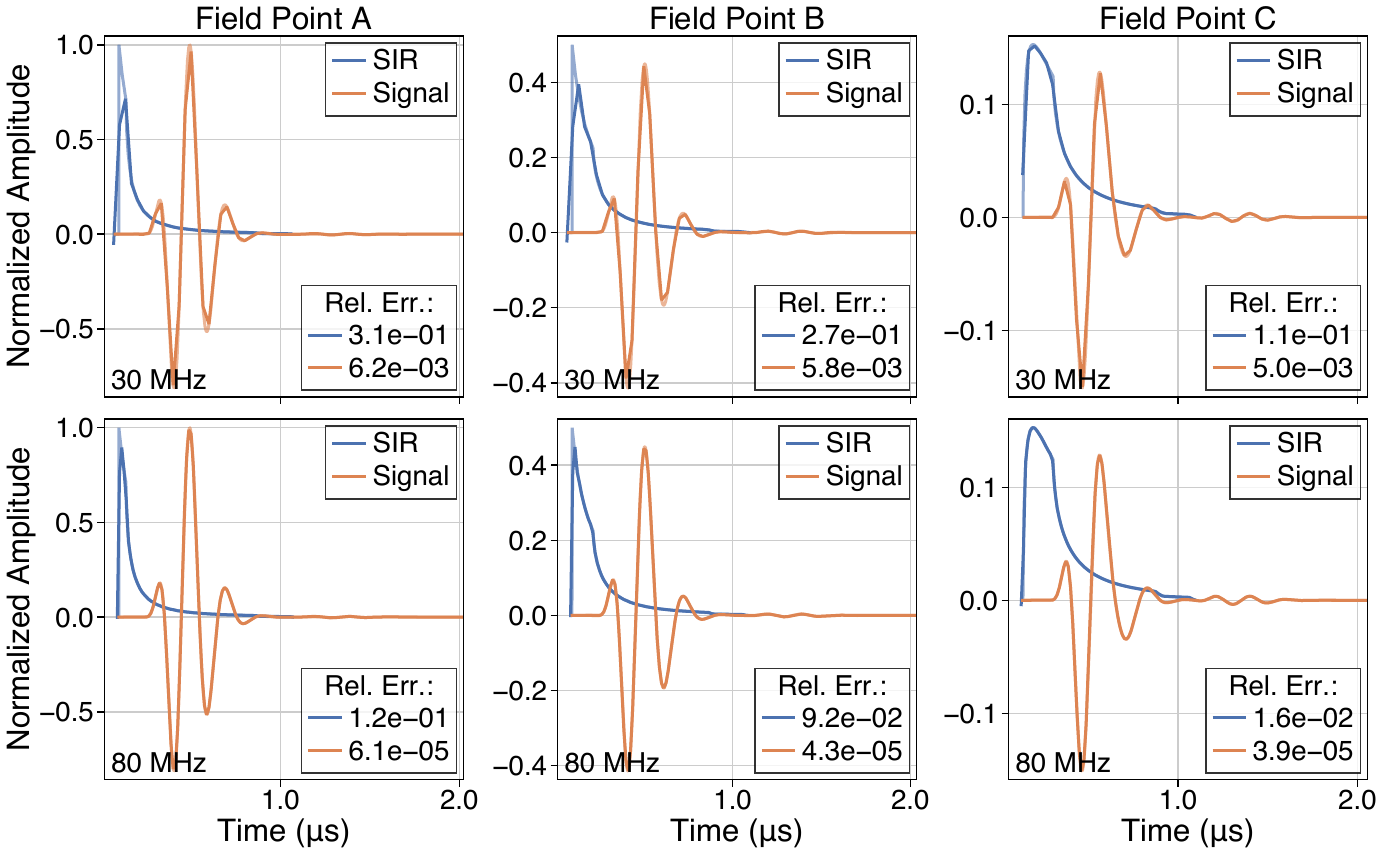}%
    \caption{(Color online) \Gls{ibid} for a cubic \gls{bspline} basis function.}%
    \label{fig:app:validation:analytic:rectangle-soft:bspline3}
\end{figure*}

\begin{figure*}[htb]
    \centering
    \includegraphics[scale=0.925]{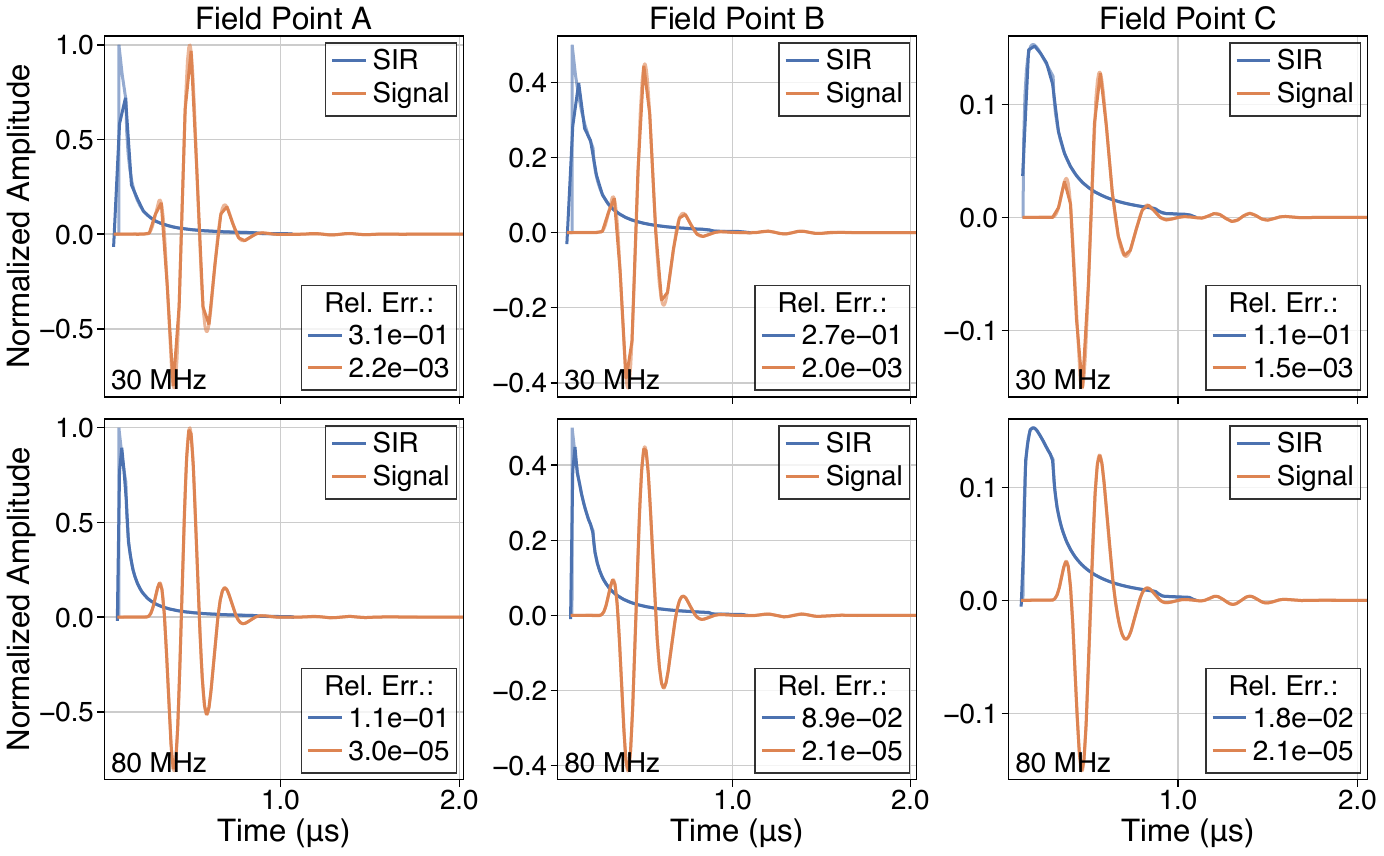}%
    \caption{(Color online) \Gls{ibid} for a cubic \glsxtrshort{o-moms} basis function.}%
    \label{fig:app:validation:analytic:rectangle-soft:omoms3}
\end{figure*}

\clearpage
\subsubsection{Rectangular element with a rigid baffle condition: Additional results}%
\label{sec:sup:validation:analytic:rectangle-rigid}

\begin{figure*}[htb]
    \centering
    \includegraphics[scale=0.925]{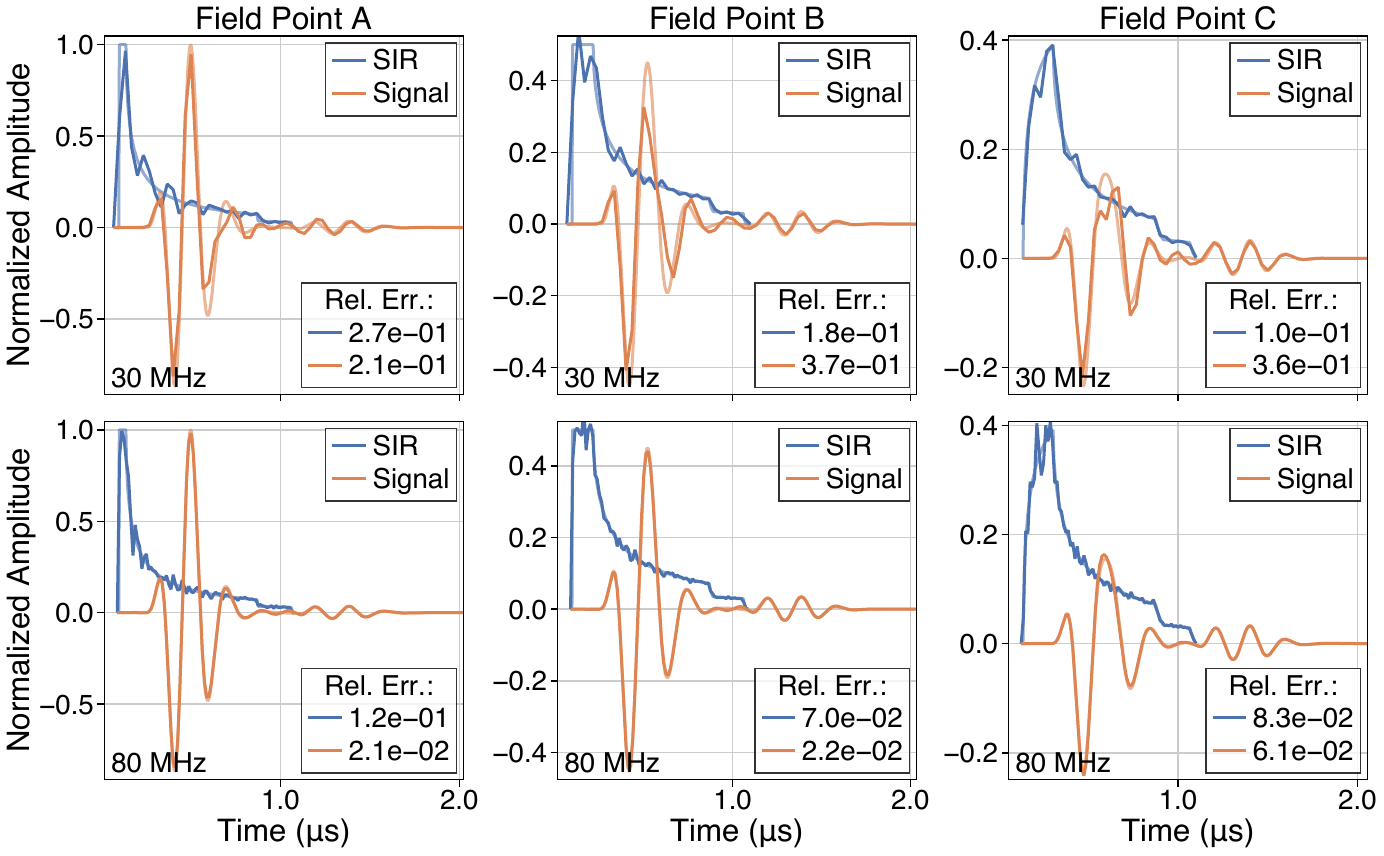}%
    \caption{%
        (Color online)
        Comparison of the \glsxtrfullpl{sir} and field signals radiated
        at different field points by
        a rectangular transducer element with a rigid baffle condition,
        excited by a windowed-sinusoidal waveform.
        The excitation waveform is a differentiated log-normal-windowed sine wave,
        with a characteristic (center) wavelength \( \lambda \).
        The geometry of the rectangular element is defined by
        a width of \( \lambda \) and a height of \( 10 \lambda \).
        The three field points (A, B, C) lie in a plane parallel to the element.
        They were positioned at a depth of \( \lambda / 2 \)
        and an elevation of \( \lambda / 2 \),
        with lateral coordinates of 0, \( \lambda / 2 \), and \( \lambda \),
        respectively.
        The proposed approach was implemented with
        a nearest-neighbor basis function (degree zero),
        and was evaluated at two sampling rates of
        (first row) \SI{30}{\mega\hertz}
        and
        (second row) \SI{80}{\mega\hertz}.
        The reference \glsxtrshortpl{sir} and field signals were evaluated at a
        sampling rate of \SI{20}{\tera\hertz}.
        They are depicted with the same colors as the approximated counterparts,
        with a lower opacity.%
    }%
    \label{fig:app:validation:analytic:rectangle-rigid:nearest}
\end{figure*}

\begin{figure*}[htb]
    \centering
    \includegraphics[scale=0.925]{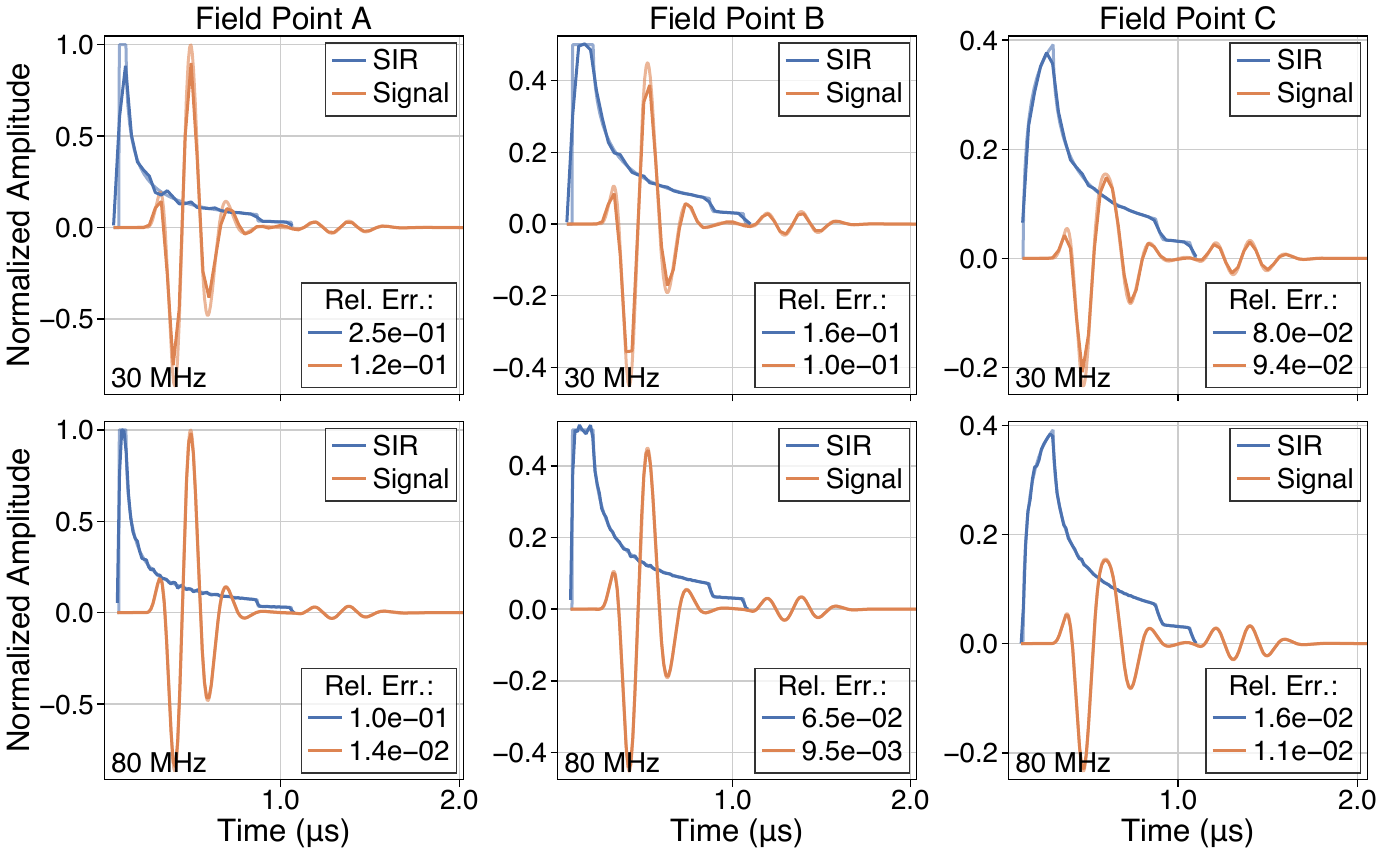}%
    \caption{(Color online) \Gls{ibid} for a linear basis function (degree one).}%
    \label{fig:app:validation:analytic:rectangle-rigid:linear}
\end{figure*}

\begin{figure*}[htb]
    \centering
    \includegraphics[scale=0.925]{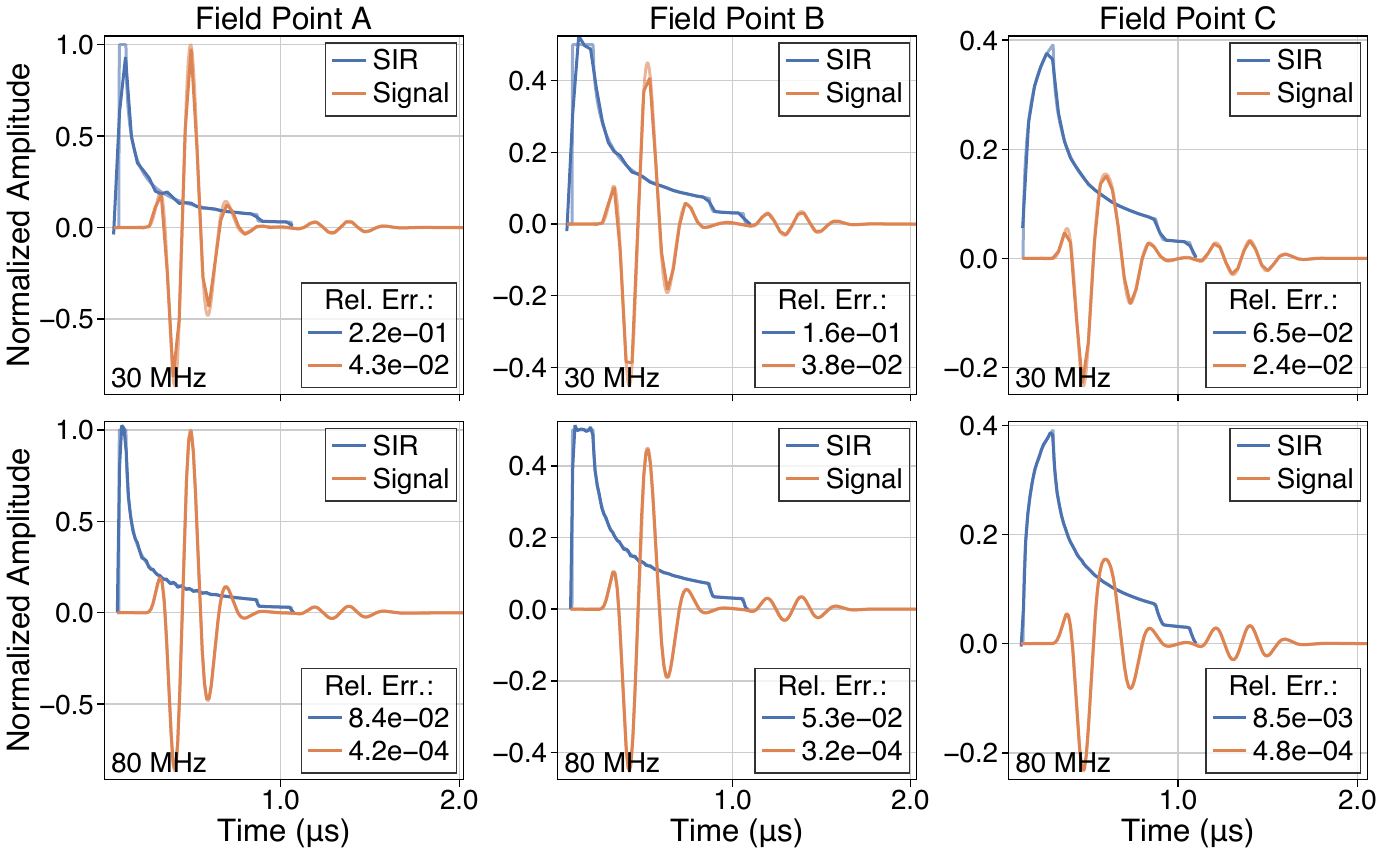}%
    \caption{(Color online) \Gls{ibid} for a quadratic \gls{keys} basis function.}%
    \label{fig:app:validation:analytic:rectangle-rigid:keys}
\end{figure*}

\begin{figure*}[htb]
    \centering
    \includegraphics[scale=0.925]{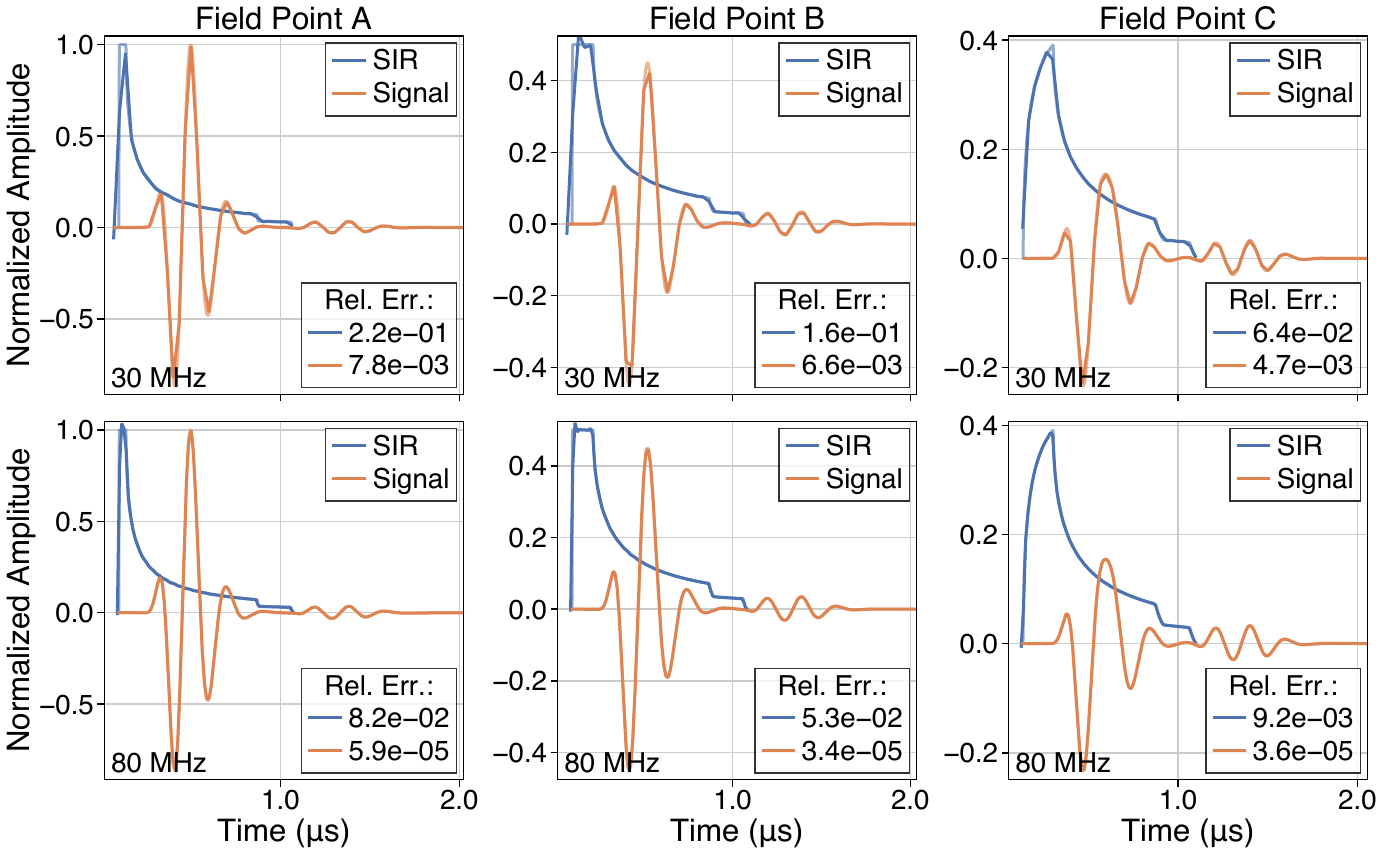}%
    \caption{(Color online) \Gls{ibid} for a cubic \gls{bspline} basis function.}%
    \label{fig:app:validation:analytic:rectangle-rigid:bspline3}
\end{figure*}

\begin{figure*}[htb]
    \centering
    \includegraphics[scale=0.925]{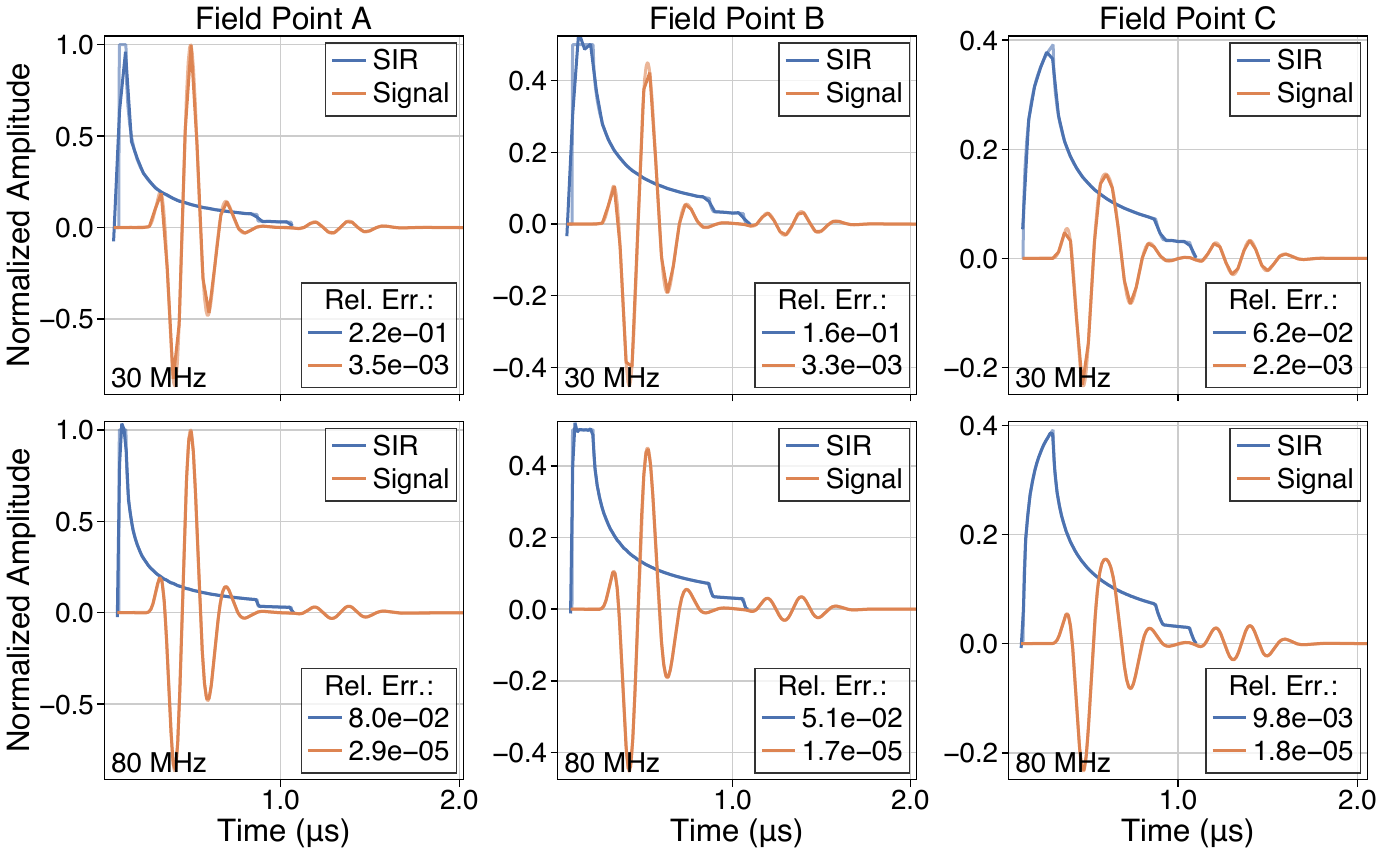}%
    \caption{(Color online) \Gls{ibid} for a cubic \glsxtrshort{o-moms} basis function.}%
    \label{fig:app:validation:analytic:rectangle-rigid:omoms3}
\end{figure*}

%% file: tables/simulation/tab-exp-rect-rigid-jasa.tex
\sisetup{
    table-format = 1.2e+1,  %
    retain-zero-exponent,
}%
\begin{tabular}{c c S S S S S S}
    \hline\hline
    {Freq.}
    & {Point}
    & {Nearest}
    & {Linear}
    & {Keys}
    & {\Gls{bspline}3}
    & {\glsxtrshort{o-moms}3}
    & {\Gls{bspline}5}
    \\
    \hline
    \multirow{3}{*}{\rotatebox[origin=c]{90}{\SI[retain-zero-exponent=false]{30}{\mega\hertz}}}
    & A
    & 2.14e-01 & 1.23e-01 & 4.34e-02 & 7.83e-03 & 3.55e-03 & 1.31e-03
    \\
    & B
    & 3.65e-01 & 1.03e-01 & 3.82e-02 & 6.58e-03 & 3.29e-03 & 1.14e-03
    \\
    & C
    & 3.59e-01 & 9.45e-02 & 2.37e-02 & 4.68e-03 & 2.24e-03 & 6.55e-04
    \\
    \hline
    \multirow{3}{*}{\rotatebox[origin=c]{90}{\SI[retain-zero-exponent=false]{80}{\mega\hertz}}}
    & A
    & 2.09e-02 & 1.39e-02 & 4.20e-04 & 5.94e-05 & 2.90e-05 & 8.28e-07
    \\
    & B
    & 2.19e-02 & 9.50e-03 & 3.19e-04 & 3.39e-05 & 1.70e-05 & 7.74e-07
    \\
    & C
    & 6.06e-02 & 1.08e-02 & 4.79e-04 & 3.58e-05 & 1.76e-05 & 3.44e-07
    \\
    \hline\hline
\end{tabular}